\documentclass{article}
\usepackage{graphicx} 

\usepackage[utf8]{inputenc}
\usepackage[english]{babel}
\usepackage{graphicx}
\usepackage{amssymb}
\usepackage{amsmath}
\usepackage{tikz}
\usetikzlibrary{decorations.pathreplacing} 
\usepackage{commath}
\usepackage{upgreek}
\usepackage[hidelinks]{hyperref}
\usepackage{color}
\usepackage{wrapfig}
\usepackage{graphicx,wrapfig,lipsum}
\graphicspath{ {./images/} }
\usepackage[a4paper, total={5.75in, 10in}]{geometry}

\usepackage{amsthm}
\newtheorem{theorem}{Theorem}[section]

\newtheorem{lemma}[theorem]{Lemma}
\newtheorem{definition}{Definition}[section]

\title{Global weak solutions to a doubly degenerate \\nutrient taxis system on the whole real line}

\author{Federico Herrero-Hervás \\
\small Universidad Complutense de Madrid, Universität Paderborn, ICAI-Universidad Pontificia Comillas}

\date{\vspace{-5ex}}

\begin{document}
\maketitle

\section*{Abstract}
This work addresses the one-dimensional Cauchy problem for the doubly degenerate nutrient taxis model
\begin{equation*}
 \begin{cases}
         \displaystyle \frac{\partial u}{\partial t} = \frac{\partial}{\partial x}(u v u_x) - \frac{\partial}{\partial x}(u^2 v v_x) + u v, & x\in \mathbb{R}, ~t>0, \\\\
       \displaystyle \frac{\partial v}{\partial t} = \frac{\partial^2 v}{\partial x^2} - u v, & x\in \mathbb{R}, ~t>0, \\\\
       u(x,0) = u_0(x) \geq 0, \quad v(x,0) = v_0(x)>0, ~  & x\in \mathbb{R},
    \end{cases}
\end{equation*}
which models pattern formation in bacterial populations. The global existence of weak solutions is established for initial data satisfying appropriate regularity and integrability conditions. To account for the degeneracy caused by $u_0$ not being strictly positive and the difficulties arising from the unboundedness of the domain, we consider a family of regularized problems posed on bounded intervals $(-\frac{1}{\varepsilon}, \frac{1}{\varepsilon})$, for $\varepsilon \in (0,1)$. Through adequate estimates uniform in $\varepsilon$, we construct global solutions to the system by passing to the limit using the Aubin-Lions lemma.

\section{Introduction}
Intricate patterns can arise in multiple scenarios in bacterial colonies as a response to changes in environmental conditions. For instance,
the availability of nutrients or the introduction of an attractant can trigger various forms of aggregation. Such patterns can be reproduced in vitro on agar plates. In this context, for certain bacterial species as \textit{Bacillus subtilis}, several studies \cite{fractal1, fractal2} have investigated the geometry of these aggregations. In \cite{ohgiwari} the shape of different aggregations is analyzed with respect to varying agar and nutrient concentrations. For hard mediums---those with high agar concentrations---in the presence of low nutrient levels, complex branching formations have been reported.
\\\\
From a mathematical point of view, to model these phenomena, nutrient taxis systems of the following form have been considered,
\begin{equation}\label{1.1}
    \begin{cases}
         \displaystyle \frac{\partial u}{\partial t} = \Delta u - \nabla \cdot (u S(u,v) \nabla v)+f(u,v), \\\\
       \displaystyle \frac{\partial v}{\partial t} = \Delta v - uv,  
    \end{cases}
\end{equation}
where $u$ represents the bacterial density and $v$ the nutrient concentration. Functions $S$ and $f$ respectively represent the taxis sensitivity coefficient and the bacterial proliferation model. System \eqref{1.1} is usually posed in bounded convex domains together with no-flux boundary conditions. This setting is essentially different from the widely studied classical Keller-Segel system 
\begin{equation}\label{1.2}
        \begin{cases}
         \displaystyle \frac{\partial u}{\partial t} = \Delta u -\chi \nabla \cdot (u \nabla v) , \\\\
       \displaystyle \frac{\partial v}{\partial t} = \Delta v +u- v,  
    \end{cases}
\end{equation}
and its extensions, in the sense that the chemotaxis migration is directed by a signal substance---in this case the nutrient---which is not produced but instead consumed by the bacteria. 
\\\\
In the case of system \eqref{1.1}, the current knowledge is much more limited, particularly concerning results capturing the above mentioned aggregation structures. In particular, for $S \equiv 1$ and $f \equiv 0$, solutions approach the spatially homogeneous steady states $u \equiv a$ and $v \equiv 0$ for appropriate $a \geq0$ \cite{tao-wink} and similar dynamics are obtained when $f = uv$ \cite{winkler-homog}. For other choices of $S$ such as $S(u,v) = \frac{1}{v}$, traveling wave solutions have been obtained in \cite{trav-waves}, but still leading to eventual spatial homogeneity.
\\\\
However, substantially different dynamics occur when bacterial diffusion is considered such that it degenerates at small signal densities. The system proposed in \cite{plaza13}, given by
\begin{equation}\label{1.3}
        \begin{cases}
         \displaystyle \frac{\partial u}{\partial t} = \nabla \cdot (uv \nabla u) - \chi \nabla \cdot (u^2 v \nabla v) + uv , \\\\
       \displaystyle \frac{\partial v}{\partial t} = \Delta v -u v,  
    \end{cases}
\end{equation}
for $\chi >0$ has recently been a source of interest. In the numerical simulations presented in the original article, the authors showed the emergence of complex patterns, which resemble the branching structures from the previously mentioned experimental records. The system was also formally derived by means of parabolic limits in \cite{plaza19}. 
\\\\
From an analytical perspective, the first work in this direction, presented in \cite{winkler1d}, successfully addressed the global solvability of system \eqref{1.3} over a bounded one-dimensional domain $\Omega$ together with no flux boundary conditions. Under the assumption that the initial data $u(x,0) = u_0(x)$ and $v(x,0) = v_0(x)$ satisfy
\begin{equation}\label{1.4}
    \begin{cases}
        u_0 \in C^\vartheta(\bar{\Omega}) \text{ for some } \vartheta \in (0,1), \text{ with } u_0 \geq 0 \text{ and } \displaystyle \int_{\Omega} \ln u_0 > -\infty,
        \\\\
        v_0 \in W^{1,\infty}(\Omega) \text{ is such that } v_0 >0 \text{ in } \bar{\Omega},
    \end{cases}
\end{equation}
globally bounded weak solutions are constructed by the limit of a sequence of regularized problems that prevent the degeneracy caused by $u_0$ not being strictly positive over $\Omega$, which is the case in most experiments. Moreover, in the same contribution, the large time behavior of the solutions is analyzed, showing that $u$ converges to a limit function in $L^\infty(\Omega)$, which is obtained as a fine time evaluation of a porous medium type equation. This behavior greatly differs from the stabilization properties that solutions of system \eqref{1.1} undergo, and thus offers an explanation to the aggregation phenomena. In \cite{li-wink}, hypotheses \eqref{1.4} are weakened, by only assuming that 
\begin{equation}\label{1.5}
    \begin{cases}
        u_0, v_0 \in  W^{1,\infty}(\Omega), \text{ with } u_0 \geq 0 \text{ and } v_0 >0 \text{ in } \bar{\Omega},
        \\\\
        \text{There exists } K >0 \text{ such that } ||u_0||_{L^\infty(\Omega)} + ||v_0||_{L^\infty(\Omega)} + ||\partial_x \ln v_0||_{L^\infty(\Omega)} \leq K.
    \end{cases}
\end{equation}
Global solvability has also been proved in convex and bounded two dimensional domains in \cite{winkler2d}, under similar hypotheses. The global boundedness of such solutions was obtained in \cite{winkler-l-inf-2d} for a broader class of systems involving \eqref{1.3}. In \cite{logistic} the system is studied with a generalized logistic growth term for the bacteria, similarly obtaining global weak solutions under an appropriate parameter relationship. 
\\\\
The aim of this work is to study the global solvability of system \eqref{1.3} posed over the whole one-dimensional space, with $\Omega = \mathbb{R}$, this is
\begin{equation}\label{sist}
    \begin{cases}
         \displaystyle \frac{\partial u}{\partial t} = \frac{\partial}{\partial x}(u v u_x) - \frac{\partial}{\partial x}(u^2 v v_x) + u v, & x\in \mathbb{R}, ~t>0, \\\\
       \displaystyle \frac{\partial v}{\partial t} = \frac{\partial^2 v}{\partial x^2} - u v, & x\in \mathbb{R}, ~t>0, \\\\
       u(x,0) = u_0(x) \geq 0, \quad v(x,0) = v_0(x)>0, ~  & x\in \mathbb{R}.
    \end{cases}
\end{equation}
For this setting, both assumptions in \eqref{1.4} and \eqref{1.5} become incompatible, as integrability properties of the initial data cannot be sustained by boundedness or integrability of their logarithms. To address the problem over the whole space, a regularized version of \eqref{sist} will be studied over balls of radius $\frac{1}{\varepsilon}$ for $\varepsilon \in (0,1)$, obtaining domain-independent estimates, with the aim of passing to the limit when $\varepsilon \searrow 0$. 
\\\\
Over the years, several other Cauchy problems have been studied for Keller-Segel related systems. For instance, in \cite{plane}, the minimal Keller-Segel system is studied over the whole two-dimensional space, obtaining a critical mass threshold for the initial data of $u$. For higher dimensions, blow up phenomena at the origin for radially-symmetric functions has been determined \cite{winkler-blow-up}, whereas boundedness can be obtained for small enough initial data in critical spaces \cite{critical1, critical2}. In \cite{preprint1} a logistic growth is incorporated, studying the global boundedness of the Cauchy problem in arbitrary dimensions. Regarding nonlinear sensitivity, in \cite{besov} a taxis term of the form $\nabla \cdot( u^m \nabla v)$ is considered, obtaining global existence in Besov spaces. However, to our knowledge, system \eqref{1.3} remains unstudied over the whole space.
\\\\
\textbf{Main results} \quad  For our case, the one-dimensional Cauchy problem \eqref{sist} will be studied under the hypotheses
\begin{equation}\label{hip0}
    \begin{cases}
        u_0 \in W^{1,\infty}(\mathbb{R}) \cap L^1(\mathbb{R}), \text{ with } u_0 \geq 0 \text{ in } \mathbb{R}, \\\\
        v_0 \in W^{1,\infty}(\mathbb{R}) \cap L^1(\mathbb{R}), \text{ with } v_0 > 0  \text{ in } \mathbb{R},
    \end{cases}
\end{equation}
as well as the additional assumptions 
\begin{equation}\label{hip01}
       \displaystyle \int_{\mathbb{R}} \left| \left(v_0^{\frac{3}{2(p+1)}}\right)_x \right|^{\frac{2(p+1)(p+2)}{p+4}} < \infty ~\text{ for all } p \geq 2,\quad \int_{\mathbb{R}} \frac{v_{0x}^2}{v_0} < \infty.
\end{equation}
Our main result guarantees that under such requirements, a global weak solution to the Cauchy problem \eqref{sist} can always be found.
\begin{theorem}\label{t1}
    Assume that \eqref{hip0} and \eqref{hip01} hold. Then problem \eqref{sist} admits a global weak solution in the sense of Definition \ref{weak-sol} below.
\end{theorem}
It is worth noting that both hypotheses in \eqref{hip01} are already included in \cite{winkler1d} in the appropriate sense over the bounded domain $\Omega$ considered. However, the boundedness of $\Omega$ ensures that \eqref{1.4} are sufficient to obtain them. In our case, the extension to the whole real line implies that both bounds have to be independently assumed, and are key for establishing Lemma \ref{l3.7} and Lemma \ref{l4.1}, respectively. Moreover, the integrability of both initial values is independently required in the whole space setting and is particularly utilized in Lemmas \ref{l3.3} to \ref{l3.5}.
\\\\
\textbf{Structure of the work} \quad After this introduction, the article begins with some preliminaries in Section \ref{s2}, which are needed in order to establish the concept of weak solutions to the problem, as well as the regularization approach followed. In particular, we consider $\varepsilon \in (0,1)$ and solutions $(u_\varepsilon, v_\varepsilon)$ to the regularized system \eqref{reg}, which allows us to deal with both degeneracies and the unboundedness of the domain. 
\\\\Next, the starting point of the analysis in Section \ref{s3} relies on a functional of the form
$$
\int_{B_{1/\varepsilon}} u_\varepsilon^p + \int_{B_{1/\varepsilon}} v_\varepsilon^{-\alpha} |v_{\varepsilon x}|^q, ~ \text{ for }\alpha = \frac{(2p-1)q}{2(p+1)} >0,
$$
for arbitrary $p \geq 2$ and  suitably chosen $q>1$. Upon suitable estimates of the quantities involved, in Lemma \ref{l3.6} we arrive at a differential inequality given by
$$
 \frac{d}{dt}  \left \{ \int_{B_{1/\varepsilon}} u_\varepsilon^p + \int_{B_{1/\varepsilon}} v_\varepsilon^{- \frac{(2p-1)q}{2(p+1)}} |v_{\varepsilon x}|^q \right \}  \leq C ||v_\varepsilon||_{L^\infty(B_{1/\varepsilon})} \Big(\int_{B_{1/\varepsilon}} u_\varepsilon^p  + 1 +\varepsilon^{\frac{p(q+2)+2}{q}} + \varepsilon^{\frac{q}{2}} \Big),
$$
which leads to bounds for $u_\varepsilon$ in $L^p(B_{1/\varepsilon})$ that are uniform in $\varepsilon$. The bounds obtained are time dependent, in contrast to those in \cite{winkler1d, li-wink}, as those rely on the assumptions on the logarithm of the initial data given in \eqref{1.4} and \eqref{1.5}, yielding the integrability of $||v_\varepsilon||_{L^\infty(B_{1/\varepsilon})}$ in time. 
\\However, the temporal dependency is enough to obtain local compactness for ultimately passing to the limit.
\\\\
The following Sections \ref{s4} and \ref{s5} are focused on obtaining compactness properties of the sequence $(u_\varepsilon^{\frac{p+1}{2}} v_\varepsilon)_{\varepsilon\in(0,1)}$, as similarly done in \cite{winkler2d}. The main difference relies on localizing the approach by means of a cutoff function as introduced in Definition \ref{ctf}, to eliminate any domain dependency. Lemma \ref{l4.5} provides a bound for $(u_\varepsilon^{\frac{p+1}{2}} v_\varepsilon)_{\varepsilon\in(0,1)}$ in $L^2\big((0,T); W^{1,1}_\text{loc}(\mathbb{R})\big)$, while Lemma \ref{l5.2} ensures that its time derivative is bounded in $L^1\big((0,T); (W^{3,2}_\text{loc}(\mathbb{R}))^*\big)$.
\\\\
Lastly, in Section \ref{s6}, an application of the Aubin-Lions lemma allows us to obtain a subsequence though with limit functions $u$ and $v$ can be defined, such that they form a weak solution to our problem \eqref{sist}. The article finishes with a proof of Theorem \ref{t1}.

\section{Preliminaries: Regularized problems and basic estimates}\label{s2}
Given the degenerate nature of system \eqref{sist}, we consider the concept of weak solutions in the following sense.
\begin{definition}\label{weak-sol}
Let $u$  and $v$ be nonnegative functions defined on $\mathbb{R} \times (0,\infty)$ satisfying
\begin{equation}\label{2.1}
    \begin{cases}
        u \in L^1_\text{loc} \left(\mathbb{R} \times [0, \infty)\right), \\\\
        v \in L^1_\text{loc} \left([0, \infty); W^{1,1}(\mathbb{R})\right),
    \end{cases}
\end{equation}
and 
\begin{equation}\label{2.2}
   u^2 v_x, ~u^2 v v_x, ~uv\in L^1_\text{loc}\left(\mathbb{R} \times [0, \infty)\right),  
\end{equation}
then, $(u,v)$ will be called a weak solution to system \eqref{sist} if 
\begin{equation}\label{2.3}
     \int_0^\infty \int_{\mathbb{R}} u \varphi_t + \int_\mathbb{R} u_0 \varphi(\cdot, 0) =  -\frac{1}{2}\int_0^\infty \int_{\mathbb{R}} u^2 v_x \varphi_x - \frac{1}{2}\int_0^\infty \int_{\mathbb{R}}u v \varphi_{xx} - \int_0^\infty \int_\mathbb{R} u^2 v v_x \varphi_x  - \int_0^\infty \int_\mathbb{R} u v \varphi,
\end{equation}
and
\begin{equation}\label{2.4}
    \int_0^\infty \int_\mathbb{R} v \varphi_t + \int_\mathbb{R} v_0 \varphi(\cdot, 0) =\int_0^\infty \int_{\mathbb{R}} v_x \varphi_x + \int_0^\infty \int_{\mathbb{R}} u v \varphi,
\end{equation}
for all $\varphi \in C_0^\infty (\mathbb{R} \times [0, \infty))$.    
\end{definition}
The main analytical challenges in system \eqref{sist} lie on the one hand, on the presence of the degenerate cross-diffusive terms, and on the other hand, on the fact that the system is posed on the unbounded domain $\mathbb{R}$. To address these issues and construct a solution, we employ a regularization approach involving a parameter $\varepsilon \in (0,1)$. Specifically, to prevent degeneracy---arising when $u_0$ is not strictly positive---we modify the initial value of $u$ by adding a small positive perturbation $\varepsilon \zeta(x)$, where $\zeta:\mathbb{R} \to \mathbb{R}^+$ is a fixed, smooth, bounded and integrable function. 
\\\\
Additionally, to deal with the unbounded spatial domain, for each  $\varepsilon \in (0,1)$ we restrict the problem to the bounded one-dimensional ball of radius $\frac{1}{\varepsilon}$, $B_{1/\varepsilon} := (-1/\varepsilon, 1/\varepsilon)$. Over each of these balls, we aim to derive estimates uniform in $\varepsilon$. 
\\\\
Thus, for any $\varepsilon \in (0,1)$, we consider the following approximate initial-boundary value problem
\begin{equation}\label{reg}
    \begin{cases}
        \displaystyle \frac{\partial u_\varepsilon}{\partial t} = \frac{\partial}{\partial x}(u_\varepsilon v_\varepsilon u_{\varepsilon x}) - \frac{\partial}{\partial x}(u_\varepsilon^2 v_\varepsilon v_{\varepsilon x}) + u_\varepsilon v_\varepsilon, & x\in B_{1/\varepsilon}, ~t>0, \\\\
       \displaystyle \frac{\partial v_\varepsilon}{\partial t} = \frac{\partial^2 v_\varepsilon}{\partial x^2} - u_\varepsilon v_\varepsilon, & x\in B_{1/\varepsilon}, ~t>0, \\\\
       u_{\varepsilon x} =  v_{\varepsilon x}  = 0, & x\in \partial B_{1/\varepsilon}, ~t>0, \\\\
       u_\varepsilon(x,0) = u_0(x) + \varepsilon \zeta(x)>0, \quad v_\varepsilon(x,0) = v_0(x)>0, ~  & x\in B_{1/\varepsilon}.
    \end{cases}
\end{equation}
Given system \eqref{reg}, for each fixed $\varepsilon \in (0,1)$ standard cross-diffusive parabolic theory can be applied to obtain local existence of solutions as follows. 
\begin{lemma}\label{l2.1}
Assume that $u_0$ and $v_0$ are such that \eqref{hip0} holds. Then for each $\varepsilon\in (0,1)$, there exists $T_{\max, \varepsilon} \in (0, \infty]$ and functions
\begin{equation*}
\begin{cases}
    u_\varepsilon \in \bigcap_{q \geq 1} C^0\left([0, T_{\max, \varepsilon}); W^{1,q}(B_{1/\varepsilon}) \right) \cap C^{2,1}(\bar{B}_{1/\varepsilon} \times \left(0, T_{\max, \varepsilon})\right), \\\\
    v_\varepsilon \in \bigcap_{q \geq 1} C^0\left([0, T_{\max, \varepsilon}); W^{1,q}(B_{1/\varepsilon}) \right) \cap C^{2,1}(\bar{B}_{1/\varepsilon} \times \left(0, T_{\max, \varepsilon})\right), 
\end{cases}
\end{equation*}
with $u_\varepsilon >0$ and $v_\varepsilon >0$ in $\bar{B}_{1/\varepsilon} \times (0,T_{\max, \varepsilon})$ such that $(u_\varepsilon, v_\varepsilon)$ solves the regularized system \eqref{reg} in the classical sense in $B_{1/\varepsilon} \times (0, T_{\max, \varepsilon})$, with the property that 
$$\text{if } T_{\max, \varepsilon} < \infty, \text{ then } \lim \sup_{t \nearrow T_{\max, \varepsilon}} ||u_\varepsilon(\cdot, t)||_{L^\infty(B_{1/\varepsilon})} = \infty.$$
\end{lemma}
\begin{proof}
    We omit the details, as the local existence of the regularized system over a bounded domain, in this case $B_{1/\varepsilon}$, is standard  following the theory developed in \cite{amann}, and has been already established in \cite{winkler1d, winkler2d, li-wink}. We refer for instance to \cite{li-wink} Lemma 2.2 for the details.
\end{proof}
Next, we prove some basic estimates for $u_\varepsilon$ and $v_\varepsilon$. 
\begin{lemma}\label{l2.2}
Assume \eqref{hip0}. Then for all $\varepsilon \in (0,1)$, $t \in (0, T_{\max, \varepsilon})$ we have
$$ \int_{B_{1/\varepsilon}} u_\varepsilon(\cdot,t) + \int_{B_{1/\varepsilon}} v_\varepsilon(\cdot,t) \leq \int_{B_{1/\varepsilon}} u_0  + \int_{B_{1/\varepsilon}} \zeta + \int_{B_{1/\varepsilon}} v_0, 
$$
as well as
$$\int_0^t \int_{B_{1/\varepsilon}} u_\varepsilon(\cdot,t) v_\varepsilon(\cdot,t) \leq \int_{B_{1/\varepsilon}} v_0,$$
and
$$||v_\varepsilon(\cdot,t)||_{L^\infty(B_{1/\varepsilon})} \leq ||v_0||_{L^\infty(B_{1/\varepsilon})} \leq ||v_0||_{L^\infty(\mathbb{R})}.$$
\end{lemma}
\begin{proof}
    Given the Neumann homogeneous boundary conditions for system \eqref{reg}, integrating both equations over $B_{1/\varepsilon}$ yields
    \begin{equation}\label{2.5}
          \frac{d}{dt} \int_{B_{1/\varepsilon}} u_\varepsilon = \int_{B_{1/\varepsilon}} u_\varepsilon v_\varepsilon, \quad \frac{d}{dt}\int_{B_{1/\varepsilon}} v_\varepsilon = - \int_{B_{1/\varepsilon}} u_\varepsilon v_\varepsilon, \quad \text{for all t } \in (0,T_{\max, \varepsilon}),  
    \end{equation}    
    so adding them and integrating in time leads to
    $$
     \int_{B_{1/\varepsilon}} u_\varepsilon(\cdot,t) + \int_{B_{1/\varepsilon}} v_\varepsilon(\cdot,t) = \int_{B_{1/\varepsilon}} u_\varepsilon(\cdot,0) + \int_{B_{1/\varepsilon}} v_\varepsilon(\cdot,0), \quad \text{for all t } \in (0,T_{\max, \varepsilon}).
    $$
    The initial value of $u_\varepsilon$ combined with the fact that $\varepsilon <1$ proves the first property. Next, a time integration of the second identity in \eqref{2.5} entails the second property, while the maximum principle applied to the second equation leads to the third one.
\end{proof}
Notice that having $L^1$ and $L^\infty$ bounds for $v_\varepsilon$ leads to $L^p$ estimates for $p \geq 2$ uniform in time and in $\varepsilon$.
\\\\
For the analysis developed in the subsequent sections, we also consider the following versions of the Gagliardo-Niremberg interpolation inequality in one dimension. In this case, the constant is domain-independent, with the dependency only present in the penalization term in a non increasing factor, if the parameters are chosen adequately. The proof strategy relies on a scaling argument, as done in \cite{kang-lee-winkler-NS20}.
\begin{lemma}\label{gn1}
    Let $p>1$, $q \in (0,p)$, $r \geq 1$, $\theta\in[0,1]$ such that
    $$\frac{1}{p} = \theta\left( \frac{1}{r}-1\right) + \frac{1-\theta}{q},$$
    and arbitrary $\sigma>0$. In dimension one, there exists $C>0$ such that for any $\varepsilon \in (0,1)$
    $$
    ||\varphi||_{L^p(B_{1/\varepsilon})} \leq C ||\varphi_x||^\theta_{L^r(B_{1/\varepsilon})} \cdot ||\varphi||_{L^q(B_{1/\varepsilon})}^{1-\theta} + C \varepsilon^{(\frac{1}{\sigma} - \frac{1}{p})} ||\varphi||_{L^\sigma(B_{1/\varepsilon})},
    $$
    for all $\varphi \in L^q(B_{1/\varepsilon})$ such that $\varphi_x \in L^r(B_{1/\varepsilon})$.
\end{lemma}
\begin{proof}
    For any given $\varphi_1 \in L^q(B_1)$ with $\varphi_x \in L^r(B_1)$, by means of the standard Gagliardo-Niremberg inequality in $B_1$, there exists $C_1 = C_1(B_1)$ such that
    \begin{equation}\label{gn1-1}
         ||\varphi_1||_{L^p(B_1)} \leq C_1 ||(\varphi_1)_x||^\theta_{L^r(B_1)} \cdot ||\varphi_1||_{L^q(B_1)}^{1-\theta} + C_1||\varphi_1||_{L^\sigma(B_1)}.
    \end{equation}
Next, to extend the inequality to $B_{1/\varepsilon}$ for an arbitrary $\varepsilon \in (0,1)$, given $\varphi \in L^q(B_{1/\varepsilon})$ with $\varphi_x \in L^r(B_{1/\varepsilon})$, we define the function
\begin{equation*}
\begin{split}
   \varphi_1:  B_1  \to & ~\mathbb{R} \\
    y \mapsto  & ~\varphi_1(y) = \varphi\left(\frac{1}{\varepsilon} \cdot y\right),
\end{split}
\end{equation*}
for which, given the one-dimensional setting, following the change of variables $x = \displaystyle \frac{1}{\varepsilon} \cdot y$ we have
$$
||\varphi_1||_{L^m(B_1)}^m = \int_{B_1} \left|\varphi\left(\frac{1}{\varepsilon} \cdot y\right)\right|^m ~dy = \varepsilon \int_{ B_{1/\varepsilon}} |\varphi(x)|^m ~dx = \varepsilon ~|| \varphi||_{L^m(B_{1/\varepsilon})}^m,
$$
for any $m>0$, as well as similarly computing the derivative $ \displaystyle ||(\varphi_1)_x||_{L^m(B_1)}^m = \varepsilon^{-(m-1)} ~||\varphi_x||_{L^m(B_{1/\varepsilon})}^m$.
\\\\
Thus, substituting in \eqref{gn1-1}, we obtain
$$
\varepsilon^{\frac{1}{p}} ~||\varphi||_{L^p(B_{1/\varepsilon})} \leq C_1 ~ \varepsilon^{- \frac{r-1}{r} \cdot \theta} ||\varphi_x||^\theta_{L^r(B_{1/\varepsilon})} \cdot \varepsilon^{\frac{1-\theta}{q}} ~||\varphi||_{L^q(B_{1/\varepsilon})}^{1-\theta} + C_1 ~\varepsilon^{\frac{1}{\sigma}} ~||\varphi||_{L^\sigma(B_{1/\varepsilon})},
$$
which upon multiplying by $\varepsilon^{-\frac{1}{p}}$ implies the result with $C = C_1(B_1)$, as 
$$ -\frac{1}{p} - \theta \cdot \left(\frac{r-1}{r} \right)+\frac{1-\theta}{q} = 0,$$
for the considered value of $\theta$.
\end{proof}
The same argument can be used in order to estimate the $L^\infty$ norm. As the proof relies on the same steps, we omit it for brevity reasons.
\begin{lemma}\label{gn2}
    Let $r \geq1$, $q >0$ and $\theta \in [0,1]$ be such that $\displaystyle \theta \left( \frac{1-r}{r}\right) + \frac{1-\theta}{q} = 0$. Then, in dimension one, there exists $C>0$ such that for any $\varepsilon\in(0,1)$
    $$
    ||\varphi||_{L^\infty(B_{1/\varepsilon})} \leq C ||\varphi_x||_{L^r(B_{1/\varepsilon})}^\theta \cdot ||\varphi||_{L^q(B_{1/\varepsilon})}^{1-\theta} + C \varepsilon^{\frac{1}{q}} ||\varphi||_{L^q(B_{1/\varepsilon})},
    $$
    for all $\varphi \in L^q(B_{1/\varepsilon})$ with $\varphi_x \in L^r(B_{1/\varepsilon})$
\end{lemma}

\section{$L^p$ estimates for $u_\varepsilon$ and global existence of regularized solutions}\label{s3}
The main objective of this section is to derive an estimate for $||u_\varepsilon||_{L^p(B_{1/\varepsilon})}$ that remains uniform with respect to $\varepsilon \in (0,1)$ for arbitrary $p \geq 2$. Although the resulting bounds will generally depend on time, they will be enough to grant the convergence in the sense of Definition \ref{weak-sol}. Moreover, these estimates will allow us to prove that regularized solutions indeed exist globally.
\\\\
To do so, we consider a functional of the following form, as in \cite{winkler1d}
\begin{equation} \label{3.1}
\int_{B_{1/\varepsilon}} u_\varepsilon^p + \int_{B_{1/\varepsilon}} v_\varepsilon^{-\alpha} |v_{\varepsilon x}|^q,
\end{equation}
with 
$$\alpha = \frac{(2p-1)q}{2(p+1)} >0,$$
for an adequately chosen $q$. We begin by computing the time derivative of \eqref{3.1}.
\begin{lemma} \label{l3.1}
    Let $p>1$, $q >2$, and $\eta > 0$. Then, there exists $C_1 > 0$ and $C_2(\eta) >0$ such that for all $\varepsilon \in (0,1)$ and $t \in (0, T_{\max, \varepsilon})$
    \begin{equation} \label{3.2}
            \frac{d}{dt} \int_{B_{1/\varepsilon}} v_\varepsilon^{-\frac{(2p-1)q}{2(p+1)}} |v_{\varepsilon x}|^q + \frac{1}{C_1} \int_{B_{1/\varepsilon}} v_\varepsilon ^{- \frac{(2p-1)q}{2(p+1)} -2} |v_{\varepsilon x}|^{q+2} \leq C_1 \int_{B_{1/\varepsilon}} u_\varepsilon ^{\frac{q+2}{2}} v_\varepsilon^{q - \frac{(2p-1)q}{2(p+1)}}, 
    \end{equation}
    and
        \begin{equation}\label{3.3}
        \begin{split}
            & \frac{d}{dt} \int_{B_{1/\varepsilon}} u_\varepsilon^p + \frac{p(p-1)}{2} \int_{B_{1/\varepsilon}} u_\varepsilon^{p-1} v_\varepsilon u_{\varepsilon x}^2 + \frac{p(p-1)}{2} \int_{B_{1/\varepsilon}} u_\varepsilon^{p+1} v_\varepsilon v_{\varepsilon x}^2   \\\\
           & \leq  p ||v_\varepsilon||_{L^\infty(B_{1/\varepsilon})}  \int_{B_{1/\varepsilon}} u_\varepsilon^p + \eta \int_{B_{1/\varepsilon}} v_\varepsilon^{-\frac{(2p-1)q}{2(p+1)} -2} |v_{\varepsilon x}|^{q+2} + C_2(\eta) \int_{B_{1/\varepsilon}} u_\varepsilon ^{\frac{(p+1)(q+2)}{q}} v_\varepsilon^{\frac{1}{q} \left(q+\frac{(2p-1)q}{(p+1)} +6\right)}               
        \end{split}
        \end{equation}
            
\end{lemma}
\begin{proof}

We omit the details for deriving \eqref{3.2}, as the result is directly given in \cite{winkler1d} Lemma 4.1 for $\frac{d}{dt} \int_{\Omega} v_\varepsilon^{-\alpha} |v_{\varepsilon x}|^q$, for a general $\Omega \subset \mathbb{R}$ and $\alpha \in (0, q)$, where the equation satisfied by $v_\varepsilon$ is the same one as here. Ours is just the particular case for $\alpha = \frac{(2p-1)q}{2(p+1)}$.
\\\\
With respect to \eqref{3.3}, integrating by parts the first equation in \eqref{reg} followed by applying Young's inequality yields for all $ t \in(0,T_{\max, \varepsilon})$
\begin{equation} \label{3.4}
    \begin{split}
        \frac{1}{p} \frac{d}{dt} \int_{B_{1/\varepsilon}} u_\varepsilon^p & + (p-1) \int_{B_{1/\varepsilon}} u_{\varepsilon}^{p-1} v_\varepsilon u_{\varepsilon x}^2 = (p-1) \int_{B_{1/\varepsilon}} u_\varepsilon^p v_\varepsilon u_{\varepsilon x} v_{\varepsilon x} +  \int_{B_{1/\varepsilon}} u_\varepsilon^p v_\varepsilon \\\\
        & \leq \frac{p-1}{2} \int_{B_{1/\varepsilon}} u_\varepsilon^{p-1} v_\varepsilon u_{\varepsilon x}^2 + \frac{p-1}{2} \int_{B_{1/\varepsilon}} u_\varepsilon^{p+1} v_\varepsilon v_{\varepsilon x}^2 +  ||v_\varepsilon||_{L^\infty(B_{1/\varepsilon})} \int_{B_{1/\varepsilon}} u_\varepsilon^p,
    \end{split}
\end{equation}
which, upon adding $ \frac{p-1}{2} \int_{B_{1/\varepsilon}} u_\varepsilon^{p+1} v_\varepsilon v_{\varepsilon x}^2$ on both sides and rearranging the terms, results in
\begin{equation}\label{3.4.2}
    \begin{split}
        \frac{1}{p} \frac{d}{dt} & \int_{B_{1/\varepsilon}} u_\varepsilon^p  + \frac{p-1}{2} \int_{B_{1/\varepsilon}} u_{\varepsilon}^{p-1} v_\varepsilon u_{\varepsilon x}^2 + \frac{p-1}{2} \int_{B_{1/\varepsilon}} u_\varepsilon^{p+1} v_\varepsilon v_{\varepsilon x}^2 
\\\\ & \leq (p-1)\int_{B_{1/\varepsilon}} u_\varepsilon^{p+1} v_\varepsilon v_{\varepsilon x}^2 +  ||v_\varepsilon||_{L^\infty(B_{1/\varepsilon})} \int_{B_{1/\varepsilon}} u_\varepsilon^p, \quad \text{for all } t \in (0,T_{\max, \varepsilon}).
    \end{split}
\end{equation}
Lastly, given $\eta >0$, again by Young's inequality with exponents $\frac{q+2}{2}$, $\frac{q+2}{q} > 1$, there exists $C_2(\eta) >0$ such that
\begin{equation} \label{3.5}
    \begin{split}
        (p-1)\int_{B_{1/\varepsilon}} u_\varepsilon^{p+1} v_\varepsilon v_{\varepsilon x}^2 = (p-1) \int_{B_{1/\varepsilon}} \left( v_\varepsilon^{-\frac{(2p-1)q}{2(p+1)} -2} |v_{\varepsilon x}|^{q+2} \right)^{\frac{2}{q+2}} u_\varepsilon^{p+1} v_\varepsilon^{\frac{1}{q+2} \left( q + 2\frac{(2p-1)q}{2(p+1)} + 6\right)} \\\\
         \leq \frac{\eta}{p} \int_{B_{1/\varepsilon}} v_\varepsilon^{\frac{(2p-1)q}{2(p+1)} - 2} |v_{\varepsilon x}|^{q+2} + C_1(\eta) \int_{B_{1/\varepsilon}} u_\varepsilon^{\frac{(p+1)(q+2)}{q}} v_\varepsilon^{\frac{1}{q} \left( q + \frac{(2p-1)q}{(p+1)} + 6\right)}.
    \end{split}
\end{equation}
Directly combining \eqref{3.4.2} and \eqref{3.5} provides \eqref{3.3} for all $t \in (0, T_{\max, \varepsilon})$.
\end{proof}
Next, using the domain-independent versions of the Gagliardo-Niremberg inequality introduced in Lemmas \ref{gn1} and \ref{gn2}, we prove the following auxiliary result.
\begin{lemma}\label{l3.2}
    Let $p>0$, then there exists $C>0$ such that for any $\varepsilon \in (0,1)$
    \begin{equation*}
        \begin{split}
    || \phi \psi^{\frac{3}{p+1}}||^{p+2}_{L^\infty(B_{1/\varepsilon})} & \leq C ||\phi||_{L^1(B_{1/\varepsilon})} ||\psi||_{L^\infty(B_{1/\varepsilon})}^{\frac{3}{p+1}} \Bigg(  ||\psi||_{L^\infty(B_{1/\varepsilon})}^2 \int_{B_{1/\varepsilon}} \phi^{p-1} \psi ~\phi_x^2 +\int_{B_{1/\varepsilon}} \phi^{p+1} \psi ~\psi_x^2 
    \\
    & + \varepsilon^{p+2}  ||\phi||_{L^1(B_{1/\varepsilon})}^{p+1} ||\psi||_{L^\infty(B_{1/\varepsilon})}^3 \Bigg)
        \end{split}
    \end{equation*}
for any positive functions $\phi, ~\psi \in C^1(\bar{B}_{1/\varepsilon})$.
\end{lemma}
\begin{proof}
    Firstly, by Lemma \ref{gn2} there exists $c_1$ independent of $\varepsilon$ satisfying 
    \begin{equation} \label{3.6}
    \begin{split}
        & \displaystyle || \phi \psi^{\frac{3}{p+1}}||^{p+2}_{L^\infty(B_{1/\varepsilon})}  = || \phi^{\frac{p+1}{2}} \psi^{\frac{3}{2}}||_{L^\infty (B_{1/\varepsilon})}^\frac{2(p+2)}{p+1}      \\\\
       & \displaystyle \quad \leq c_1 ||(\phi^{\frac{p+1}{2}} \psi^{\frac{3}{2}})_x||_{L^2(B_{1/\varepsilon})}^2 \cdot ||\phi^{\frac{p+1}{2}} \psi^{\frac{3}{2}}||_{L^\frac{2}{p+1} (B_{1/\varepsilon})}^\frac{2}{p+1} +  c_1 \varepsilon^{p+2}~||\phi^{\frac{p+1}{2}} \psi^{\frac{3}{2}}||_{L^{\frac{2}{p+1}}(B_{1/\varepsilon})}^{\frac{2(p+2)}{p+1}}.
    \end{split}
    \end{equation}  
    Computing the derivative appearing in the first term on the right hand side we have
    \begin{equation*}
        \begin{split}
        ||(\phi^{\frac{p+1}{2}} \psi^{\frac{3}{2}})_x||_{L^2(B_{1/\varepsilon})}^2  &= \int_{B_{1/\varepsilon}}  \left( \frac{p+1}{2} \phi^{\frac{p+1}{2}} \psi^{\frac{3}{2}} \phi_x + \frac{3}{2} \phi^{\frac{p+1}{2}} \psi^{\frac{1}{2}} \psi_x \right)^2 \\\\
        & \leq \frac{(p+1)^2}{2} \int_{B_{1/\varepsilon}} \phi^{p-1} \psi^3 \phi_x^2 + \frac{9}{2} \int_{B_{1/\varepsilon}} \phi^{p+1} \psi ~\psi_x^2 \\\\
        & \leq \frac{(p+1)^2}{2} ||\psi||_{L^\infty(B_{1/\varepsilon})}^2 \int_{B_{1/\varepsilon}} \phi^{p-1} \psi ~\phi_x^2 + \frac{9}{2} \int_{B_{1/\varepsilon}} \phi^{p+1} \psi ~\psi_x^2,
        \end{split} 
    \end{equation*}
    and moreover for the other two terms appearing in \eqref{3.6}
    $$
    ||\phi^{\frac{p+1}{2}} \psi^{\frac{3}{2}}||_{L^{\frac{2}{p+1}}(B_{1/\varepsilon})}^{\frac{2}{p+1}} = \int_{B_{1/\varepsilon}} \phi ~\psi^{\frac{3}{p+1}} \leq ||\phi||_{L^1(B_{1/\varepsilon})} ||\psi||_{L^\infty(B_{1/\varepsilon})}^{\frac{3}{p+1}}.
    $$
    Substituting these bounds back in \eqref{3.6} leads to
    \begin{equation*}
        \begin{split}
            || \phi \psi^{\frac{3}{p+1}}||^{p+2}_{L^\infty(B_{1/\varepsilon})}  & \leq
            c_1 ||\phi||_{L^1(B_{1/\varepsilon})} ||\psi||_{L^\infty(B_{1/\varepsilon})}^{\frac{3}{p+1}} \left( \frac{(p+1)^2}{2} ||\psi||_{L^\infty(B_{1/\varepsilon})}^2 \int_{B_{1/\varepsilon}} \phi^{p-1} \psi ~\phi_x^2 + \frac{9}{2} \int_{B_{1/\varepsilon}} \phi^{p+1} \psi ~\psi_x^2  \right)
            \\\\ & + c_1 \varepsilon^{p+2}  ||\phi||_{L^1(B_{1/\varepsilon})}^{p+2} ||\psi||_{L^\infty(B_{1/\varepsilon})}^{\frac{3(p+2)}{p+1}},
        \end{split}
    \end{equation*}
    which proves the lemma, as $\displaystyle \frac{3(p+2)}{p+1} = \frac{3}{p+1} +3$.
\end{proof}
Next, using Lemma \ref{l3.2} we prove a technical result that will later be used to control the right hand side terms in \eqref{3.2} and \eqref{3.3} in order to estimate the time derivative of the functional defined in \eqref{3.1}.
\begin{lemma} \label{l3.3}
    Let $p,~r>0$ be such that 
    $$\frac{(p+1)(p+2)}{p+4} \leq r < p+2.$$
    Then, for any $\eta>0$ and $K>0$ there exists $C(\eta,K)>0$ such that whenever \eqref{hip0} holds as well as
    \begin{equation} \label{hip1}
        \int_\mathbb{R} u_0 \leq K, \quad  \int_\mathbb{R} v_0 \leq K, \quad ||v_0||_{L^\infty(\mathbb{R} )} \leq K, \text{ and} \quad \int_{\mathbb{R}} \zeta \leq K,
    \end{equation}
    we have that for any $\varepsilon \in (0,1)$, $t \in (0, T_{\max, \varepsilon})$
    $$
   ||u_\varepsilon v_\varepsilon^{\frac{3}{p+1}}||^{r}_{L^\infty(B_{1/\varepsilon})} \leq \eta \left( \int_{B_{1/\varepsilon}} u_\varepsilon^{p-1} v_\varepsilon ~u_{\varepsilon x}^2 +\int_{B_{1/\varepsilon}} u_\varepsilon^{p+1} v_\varepsilon ~ v_{\varepsilon x}^2 \right) + C(\eta,K)\Big(1+\varepsilon^r \Big) || v_\varepsilon||_{L^\infty(B_{1/\varepsilon})}
    $$
\end{lemma}
\begin{proof}
Combining the estimates proved in Lemma \ref{l2.2} with assumption \eqref{hip1}, yields for all $t \in (0, T_{\max, \varepsilon})$ 
$$
||u_\varepsilon||_{L^1(B_{1/\varepsilon})} \leq \int_{B_{1/\varepsilon}} u_0 + \int_{B_{1/\varepsilon}} \zeta  + \int_{B_{1/\varepsilon}} v_0 \leq 3K, \quad ||v_\varepsilon(\cdot,t)||_{L^\infty(B_{1/\varepsilon})} \leq ||v_0||_{L^\infty(B_{1/\varepsilon})} \leq K.
$$
Next, applying Lemma \ref{l3.2} to $\phi = u_\varepsilon, ~ \psi = v_\varepsilon$, we obtain that for any $t \in (0, T_{\max, \varepsilon})$
\begin{equation}\label{3.8}
    \begin{split}
  || u_\varepsilon v_\varepsilon^{\frac{3}{p+1}}||^{p+2}_{L^\infty(B_{1/\varepsilon})} & \leq C \cdot 3K \cdot ||v_\varepsilon||_{L^\infty(B_{1/\varepsilon})}^{\frac{3}{p+1}} \Bigg(  K^2 \int_{B_{1/\varepsilon}} u_\varepsilon^{p-1} v_\varepsilon ~u_{\varepsilon x}^2 +\int_{B_{1/\varepsilon}} u_\varepsilon^{p+1} v_\varepsilon ~ v_{\varepsilon x}^2 
    \\
    & + \varepsilon^{p+2}  (2K)^{p+1} ||v_\varepsilon||_{L^\infty(B_{1/\varepsilon})}^3 \Bigg)
    \\
    & \leq c_1 ||v_\varepsilon||_{L^\infty(B_{1/\varepsilon})}^{\frac{3}{p+1}} I(t) + c_1  \varepsilon^{p+2}  ||v_\varepsilon||_{L^\infty(B_{1/\varepsilon})}^{\frac{3(p+2)}{p+1}},
    \end{split}
\end{equation}
for a certain $c_1 = c_1(K)>0$, where 
$$
I(t) := \int_{B_{1/\varepsilon}} u_\varepsilon^{p-1} v_\varepsilon ~u_{\varepsilon x}^2 +\int_{B_{1/\varepsilon}} u_\varepsilon^{p+1} v_\varepsilon ~ v_{\varepsilon x}^2 .
$$
Thus, for any $r$ satisfying $\frac{(p+1)(p+2)}{p+4} \leq r < p+2$,  we can bound
\begin{equation}\label{3.9}
     || u_\varepsilon v_\varepsilon^{\frac{3}{p+1}}||^{r}_{L^\infty(B_{1/\varepsilon})} \leq c_2  || v_\varepsilon||^{\frac{3r}{(p+1)(p+2)}}_{L^\infty(B_{1/\varepsilon})} I(t)^{\frac{r}{p+2}} + c_2 \varepsilon^r  ||v_\varepsilon||_{L^\infty(B_{1/\varepsilon})}^{\frac{3r}{p+1}},
\end{equation}
for a certain $c_2>0$. The assumptions on $r$ grants that
$$ 
b := \frac{3r}{p+1}  -1\geq \frac{2(p+1)}{p+4} > 0,
$$
with which for all $t\in (0, T_{\max, \varepsilon})$ we obtain
$$
\varepsilon^r ||v_\varepsilon||_{L^\infty(B_{1/\varepsilon})}^{\frac{3r}{p+1}} = \varepsilon^r  ||v_\varepsilon||_{L^\infty(B_{1/\varepsilon})}^b ||v_\varepsilon||_{L^\infty(B_{1/\varepsilon})} \leq \varepsilon^r K^b ||v_\varepsilon||_{L^\infty(B_{1/\varepsilon})}.
$$
The remainder of the proof follows the same steps as Lemma 4.4 in \cite{winkler1d}. By Young's inequality, for any $\eta>0$ there exists $c_3 = c_3(\eta, K) >0$ such that
$$
c_2  || v_\varepsilon||^{\frac{3r}{(p+1)(p+2)}}_{L^\infty(B_{1/\varepsilon})} I(t)^{\frac{r}{p+2}} \leq \eta I(t) + c_3 || v_\varepsilon||^{\frac{3r}{(p+1)(p+2-r)}}_{L^\infty(B_{1/\varepsilon})} \leq \eta I(t) + c_3 K^a || v_\varepsilon||_{L^\infty(B_{1/\varepsilon})},
$$
where again $ a := \frac{3r}{(p+1)(p+2-r)} -1 $ is nonnegative by the choice of $r$. Direct substitution into \eqref{3.9} gives the result.
\end{proof}
Next, the estimate provided by Lemma \ref{l3.3} can be used to bound the right hand side term in \eqref{3.2}.
\begin{lemma} \label{l3.4}
    Let $\displaystyle p > \frac{1}{2}$, $q \in \big(1, 2(p+2)\big)$ and $\displaystyle \alpha = \frac{(2p-1)q}{2(p+1)}>0$ be such that
    \begin{equation} \label{3.10}
        \alpha \leq q - \frac{3(p+2)}{p+4}, \quad \alpha > q - \frac{3(p+2)}{p+1}.
    \end{equation}
    Then, for all $\eta >0$ and $K>0$, there exists $C(\eta, K)>0$ such that if \eqref{hip0} and \eqref{hip1} hold, the following inequality is satisfied for all $\varepsilon \in (0,1)$ and all $t \in (0, T_{\max, \varepsilon})$
    $$
    \int_{B_{1/\varepsilon}} u_\varepsilon^{\frac{q+2}{2}} v_\varepsilon^{q- \alpha}  \leq \eta \left( \int_{B_{1/\varepsilon}} u_\varepsilon^{p-1} v_\varepsilon ~u_{\varepsilon x}^2 +\int_{B_{1/\varepsilon}} u_\varepsilon^{p+1} v_\varepsilon ~ v_{\varepsilon x}^2 \right) + C(\eta,K)\Big(1+\varepsilon^{\frac{(p+1)(q-\alpha)}{3}} \Big) || v_\varepsilon||_{L^\infty(B_{1/\varepsilon})}.
    $$
\end{lemma}
\begin{proof}
    The result follows directly by noting that for $t \in (0, T_{\max, \varepsilon})$
    \begin{equation}\label{3.11}    
    \begin{split}
            \int_{B_{1/\varepsilon}} u_\varepsilon^{\frac{q+2}{2}} v_\varepsilon^{q- \alpha} = \int_{B_{1/\varepsilon}} \left(u_\varepsilon ~v_\varepsilon^{\frac{3}{p+1}} \right)^\frac{(p+1)(q-\alpha)}{3} \cdot ~ u_\varepsilon^{\frac{q+2}{2} - \frac{(p+1)(q-\alpha)}{3}} \\\\
            \leq ||u_\varepsilon ~v_\varepsilon^{\frac{3}{p+1}}||_{L^\infty(B_{1/\varepsilon})}^{\frac{(p+1)(q-\alpha)}{3}} \cdot \int_{B_{1/\varepsilon}} u_\varepsilon^{\frac{q+2}{2} - \frac{(p+1)(q-\alpha)}{3}},
    \end{split}
    \end{equation}
    where precisely, given our choice of $\alpha$, the exponent appearing in the last term is
    $$
    \frac{q+2}{2} - \frac{(p+1)(q-\alpha)}{3} = \frac{q+2}{2} - \frac{(p+1)}{3} \cdot \left(q-\frac{(2p-1)q}{2(p+1)} \right)  = \frac{q+2}{2} - \frac{q}{2} = 1,
    $$
    and therefore 
    \begin{equation}\label{3.11.2}  
    \int_{B_{1/\varepsilon}} u_\varepsilon^{\frac{q+2}{2} - \frac{(p+1)(q-\alpha)}{3}} =  ||u_\varepsilon ||_{L^1(B_{1/\varepsilon})} \leq 3K,           
    \end{equation}
    by Lemma \ref{l2.2} and hypothesis \eqref{hip1}.
    \\\\
    Moreover, assumptions \eqref{3.10} imply that
    $$
    \frac{(p+1)(q-\alpha)}{3} \geq \frac{p+1}{3} \cdot \frac{3(p+2)}{p+4} = \frac{(p+1)(p+4)}{p+4},
    $$
    as well as
    $$
    \frac{(p+1)(q-\alpha)}{3} < \frac{p+1}{3} \cdot \frac{3(p+2)}{p+1} = p+2.
    $$
    In this way, by Lemma \ref{l3.3}, considering $r := \frac{(p+1)(q-\alpha)}{3}$, for any choice of $\eta$ and $K>0$, there exists $C_1(\eta,K)>0$ such that
    $$
    ||u_\varepsilon v_\varepsilon^{\frac{3}{p+1}}||^{ \frac{(p+1)(q-\alpha)}{3}}_{L^\infty(B_{1/\varepsilon})} \leq \frac{\eta}{3K} \left( \int_{B_{1/\varepsilon}} u_\varepsilon^{p-1} v_\varepsilon ~u_{\varepsilon x}^2 +\int_{B_{1/\varepsilon}} u_\varepsilon^{p+1} v_\varepsilon ~ v_{\varepsilon x}^2 \right) + C_1(\eta,K)\Big(1+\varepsilon^{ \frac{(p+1)(q-\alpha)}{3}} \Big) || v_\varepsilon||_{L^\infty(B_{1/\varepsilon})}.
    $$
    Combining this with \eqref{3.11.2} and substituting into \eqref{3.11} finishes the proof by defining $C(\eta,K) := 3K \cdot C_1(\eta,K)$.
\end{proof}
Lastly, a similar argument can be followed for estimating the last term in \eqref{3.3}.
\begin{lemma} \label{l3.5}
Let $\displaystyle p > \frac{1}{2}$, $q >1$ and $\displaystyle \alpha = \frac{(2p-1)q}{2(p+1)}>0$ be such that
    \begin{equation} \label{3.12}
        \alpha \geq \frac{3(p+2)q}{2(p+4)} - \frac{q+6}{2},  \quad \alpha < \frac{3(p+2)q}{2(p+1)} - \frac{q+6}{2}.
    \end{equation}
Then, for all $\eta >0$ and $K>0$, there exists $C(\eta, K)>0$ such that whenever \eqref{hip0} and \eqref{hip1} are satisfied, the following estimate holds for all $\varepsilon \in (0,1)$ and all $t \in (0, T_{\max, \varepsilon})$
    \begin{equation*}
        \begin{split}
           \int_{B_{1/\varepsilon}} u_\varepsilon^{\frac{(p+1)(q+2)}{q}} v_\varepsilon^{\frac{q+2\alpha+6}{q}}  \leq \eta & \left( \int_{B_{1/\varepsilon}} u_\varepsilon^{p-1} v_\varepsilon ~u_{\varepsilon x}^2 +\int_{B_{1/\varepsilon}} u_\varepsilon^{p+1} v_\varepsilon ~ v_{\varepsilon x}^2 \right) 
           \\&+ C(\eta,K)\Big(1+\varepsilon^{\frac{(p+1)(q+2\alpha +6)}{3q}} \Big) || v_\varepsilon||_{L^\infty(B_{1/\varepsilon})}.  
        \end{split}
    \end{equation*} 
\end{lemma}
\begin{proof}
The proof follows the same steps as in Lemma \ref{l3.4}. Firstly, for all $t\in(0, T_{\max, \varepsilon})$ we have
\begin{equation}\label{3.13}
    \begin{split} 
        \int_{B_{1/\varepsilon}} u_\varepsilon^{\frac{(p+1)(q+2)}{q}} v_\varepsilon^{\frac{q+2\alpha+6}{q}}  \leq \int_{B_{1/\varepsilon}} \left(u_\varepsilon ~v_\varepsilon^{\frac{3}{p+1}} \right)^\frac{(p+1)(q+2\alpha +6)}{3q} \cdot ~ u_\varepsilon^{\frac{2(p+1)(q-\alpha)}{3q}} \\
         \leq ||u_\varepsilon ~v_\varepsilon^{\frac{3}{p+1}}||_{L^\infty(B_{1/\varepsilon})}^\frac{(p+1)(q+2\alpha +6)}{3q} \cdot \int_{B_{1/\varepsilon}} u_\varepsilon^{\frac{2(p+1)(q-\alpha)}{3q}},         
    \end{split}
\end{equation}
where again, for the chosen $\alpha$
$$
\frac{2(p+1)(q-\alpha)}{3q} = \frac{2(p+1)}{3q} \left(q- \frac{(2p-1)q}{2(p+1)}  \right) = \frac{2(p+1)}{3q} \cdot \frac{3q}{2(p+1)} = 1,
$$
therefore having a global bound for the last term in \eqref{3.13}. Lastly, the conditions on $p$, $q$ and $\alpha$ in \eqref{3.12} ensure that Lemma \ref{l3.3} can be applied again to $||u_\varepsilon ~v_\varepsilon^{\frac{3}{p+1}}||_{L^\infty(B_{1/\varepsilon})}^\frac{(p+1)(q+2\alpha +6)}{3q}$, concluding the proof.
\end{proof}
A combination of these lemmas allows us to estimate the time derivative of the functional \eqref{3.1} for the selected $\alpha$. This eventually leads us to the desired $L^p$ bounds for $u_\varepsilon$. It is important however to select a range of values of $p$ and $q$ such that assumptions \eqref{3.10} and \eqref{3.12} are simultaneously met.
\begin{lemma}\label{l3.6}
    Let $p \geq 2$ and $q\geq 4$ be such that
    $$\frac{2(p+1)(p+2)}{p+4} \leq q < 2(p+2).$$
    Then for all $K>0$ there exists $C(K)>0$ such that if \eqref{hip0} and \eqref{hip1} are satisfied, we have that for all $\varepsilon \in (0,1)$ and all $t \in (0, T_{\max, \varepsilon})$
    $$
        \frac{d}{dt}  \left \{ \int_{B_{1/\varepsilon}} u_\varepsilon^p + \int_{B_{1/\varepsilon}} v_\varepsilon^{- \frac{(2p-1)q}{2(p+1)}} |v_{\varepsilon x}|^q \right \}  \leq C(K) ||v_\varepsilon||_{L^\infty(B_{1/\varepsilon})} \Big(\int_{B_{1/\varepsilon}} u_\varepsilon^p  + 1 +\varepsilon^{\frac{p(q+2)+2}{q}} + \varepsilon^{\frac{q}{2}} \Big).
    $$
\end{lemma}
\begin{proof}
    Fixing again $\displaystyle \alpha = \frac{(2p-1)q}{2(p+1)}$, beginning with Lemma \ref{l3.1}, by \eqref{3.2} there exist $c_1, c_2 >0$ (being $c_2 = c_1^{-1}$) such that
    \begin{equation} \label{3.14}
        \frac{d}{dt} \int_{B_{1/\varepsilon}} v_\varepsilon^{-\alpha} |v_{\varepsilon x}|^q + c_1 \int_{B_{1/\varepsilon}} v_\varepsilon ^{- \alpha -2} |v_{\varepsilon x}|^{q+2} \leq c_2 \int_{B_{1/\varepsilon}} u_\varepsilon ^{\frac{q+2}{2}} v_\varepsilon^{q - \alpha}, \quad \text{for all } t \in (0,T_{\max,\varepsilon}). 
    \end{equation}
    Next, for $\eta = c_1>0$, \eqref{3.3} provides $C_2(\eta) >0$ such that for    
    $c_3 := \frac{p(p-1)}{2} >0$ and $c_4 := \max\{p ,C_2(\eta) >0\}$, we have
      \begin{equation}\label{3.15}
        \begin{split}
            & \frac{d}{dt} \int_{B_{1/\varepsilon}} u_\varepsilon^p + c_3 \left( \int_{B_{1/\varepsilon}} u_\varepsilon^{p-1} v_\varepsilon u_{\varepsilon x}^2 + \int_{B_{1/\varepsilon}} u_\varepsilon^{p+1} v_\varepsilon v_{\varepsilon x}^2 \right)  \\
           & \leq  c_4 ||v_\varepsilon||_{L^\infty(B_{1/\varepsilon})}  \int_{B_{1/\varepsilon}} u_\varepsilon^p + c_1 \int_{B_{1/\varepsilon}} v_\varepsilon^{-\alpha -2} |v_{\varepsilon x}|^{q+2} + c_4 \int_{B_{1/\varepsilon}} u_\varepsilon ^{\frac{(p+1)(q+2)}{q}} v_\varepsilon^{\frac{q+2\alpha + 6}{q}}.             
        \end{split}
    \end{equation}
    One can check that the range of values for $q$ entails \eqref{3.10}, and having $p \geq 2$ implies \eqref{3.12}. Thus, we can apply Lemmas \ref{l3.4} and \ref{l3.5}. In particular, for $\displaystyle \eta = \frac{c_3}{2 c_4}>0$, by Lemma \ref{l3.5} there exists $c_5 = c_5(K)>0$ such that for all $t \in (0, T_{\max, \varepsilon})$
    \begin{equation}\label{3.16}
    \begin{split}
       c_4 \int_{B_{1/\varepsilon}} u_\varepsilon^{\frac{(p+1)(q+2)}{q}} v_\varepsilon^{\frac{q+2\alpha+6}{q}}  \leq \frac{c_3}{2}\left( \int_{B_{1/\varepsilon}} u_\varepsilon^{p-1} v_\varepsilon ~u_{\varepsilon x}^2 +\int_{B_{1/\varepsilon}} u_\varepsilon^{p+1} v_\varepsilon ~ v_{\varepsilon x}^2 \right) 
       \\+ c_5 \Big(1+\varepsilon^{\frac{(p+1)(q+2\alpha +6)}{3q}} \Big) || v_\varepsilon||_{L^\infty(B_{1/\varepsilon})}.        
    \end{split}
    \end{equation}
    In the same way, for $\eta = \frac{c_3}{2c_2} >0$, by Lemma \ref{l3.4} we obtain $c_6 = c_6(K)>0$ for which, for all $t \in (0, T_{\max, \varepsilon})$ we have
    \begin{equation}\label{3.17} 
        c_2 \int_{B_{1/\varepsilon}} u_\varepsilon^{\frac{q+2}{2}} v_\varepsilon^{q- \alpha}  \leq \frac{c_3}{2} \left( \int_{B_{1/\varepsilon}} u_\varepsilon^{p-1} v_\varepsilon ~u_{\varepsilon x}^2 +\int_{B_{1/\varepsilon}} u_\varepsilon^{p+1} v_\varepsilon ~ v_{\varepsilon x}^2 \right) + c_6 \Big(1+\varepsilon^{\frac{(p+1)(q-\alpha)}{3}} \Big) || v_\varepsilon||_{L^\infty(B_{1/\varepsilon})}.
    \end{equation}
    Lastly, upon combining \eqref{3.14}-\eqref{3.17}, we get
    \begin{equation*}
        \begin{split}
            \frac{d}{dt} \left \{ \int_{B_{1/\varepsilon}} u_\varepsilon^p + \int_{B_{1/\varepsilon}} v_\varepsilon^{-\alpha} |v_{\varepsilon x}|^q \right \} + c_1 \int_{B_{1/\varepsilon}} v_\varepsilon ^{- \alpha -2} |v_{\varepsilon x}|^{q+2} +  c_3 \left( \int_{B_{1/\varepsilon}} u_\varepsilon^{p-1} v_\varepsilon u_{\varepsilon x}^2 + \int_{B_{1/\varepsilon}} u_\varepsilon^{p+1} v_\varepsilon v_{\varepsilon x}^2 \right) \\\\
            \leq c_2 \int_{B_{1/\varepsilon}} u_\varepsilon ^{\frac{q+2}{2}} v_\varepsilon^{q - \alpha} +  c_4 ||v_\varepsilon||_{L^\infty(B_{1/\varepsilon})}  \int_{B_{1/\varepsilon}} u_\varepsilon^p + c_1 \int_{B_{1/\varepsilon}} v_\varepsilon^{-\alpha -2} |v_{\varepsilon x}|^{q+2} + c_4 \int_{B_{1/\varepsilon}} u_\varepsilon ^{\frac{(p+1)(q+2)}{q}} v_\varepsilon^{\frac{q+2\alpha + 6}{q}}    \\\\
            \leq c_1 \int_{B_{1/\varepsilon}} v_\varepsilon^{-\alpha -2} |v_{\varepsilon x}|^{q+2}  + 2 \cdot \frac{c_3}{2} \left( \int_{B_{1/\varepsilon}} u_\varepsilon^{p-1} v_\varepsilon ~u_{\varepsilon x}^2 +\int_{B_{1/\varepsilon}} u_\varepsilon^{p+1} v_\varepsilon ~ v_{\varepsilon x}^2 \right) + c_4 ||v_\varepsilon||_{L^\infty(B_{1/\varepsilon})}  \int_{B_{1/\varepsilon}} u_\varepsilon^p \\\\
            + \left[ c_5 \Big(1+\varepsilon^{\frac{(p+1)(q+2\alpha +6)}{3q}} \Big)  + c_6 \Big(1+\varepsilon^{\frac{(p+1)(q-\alpha)}{3}} \Big) \right] || v_\varepsilon||_{L^\infty(B_{1/\varepsilon})}, \quad \text{for all } t \in (0, T_{\max, \varepsilon}).
        \end{split}
    \end{equation*}
    Equivalently, for $t \in (0, T_{\max, \varepsilon})$ we obtain
    \begin{equation}\label{3.18}
            \frac{d}{dt}  \left \{ \int_{B_{1/\varepsilon}} u_\varepsilon^p + \int_{B_{1/\varepsilon}} v_\varepsilon^{-\alpha} |v_{\varepsilon x}|^q \right \}\leq c_7 ||v_\varepsilon||_{L^\infty(B_{1/\varepsilon})} \Big( \int_{B_{1/\varepsilon}} u_\varepsilon^p + 1 +\varepsilon^{\frac{(p+1)(q+2\alpha +6)}{3q}} + \varepsilon^{\frac{(p+1)(q-\alpha)}{3}} \Big),
    \end{equation}
    and lastly, substituting $\alpha$ by its value, $\displaystyle \frac{(2p-1)q}{2(p+1)}$, yields the result.     
\end{proof}
Lastly, as a conclusion to the previous lemmas, we obtain a time-dependent $L^p$ bound for $u_\varepsilon$ for arbitrary $p \geq 2$.
\begin{lemma}\label{l3.7}
    Let $p \geq 2$, $K >0$ and assume $u_0$ and $v_0$ satisfy \eqref{hip0} and \eqref{hip1} as well as 
    \begin{equation} \label{3.19}
        \int_{\mathbb{R}} u_0^p \leq K, \quad \int_{\mathbb{R}} \left| \left(v_0^{\frac{3}{2(p+1)}}\right)_x \right|^{\frac{2(p+1)(p+2)}{p+4}} \leq K.
    \end{equation}
    Then, for all $\varepsilon \in (0,1)$, $t  \in (0,T_{\max, \varepsilon})$ there exists $C(p,k,t)>0$ independent of $\varepsilon$ such that
    $$
    \int_{B_{1/\varepsilon}} u_\varepsilon^p(\cdot, t) \leq C(p,K,t).
    $$
    
\end{lemma}
\begin{proof}
    By Lemma \ref{l3.6}, with $q = \displaystyle \frac{2(p+1)(p+2)}{p+4} \geq 4$, there exists $C(K)>0$ such that for all $t \in (0, T_{\max, \varepsilon})$
    $$
    \frac{d}{dt}  \left \{ \int_{B_{1/\varepsilon}} u_\varepsilon^p + \int_{B_{1/\varepsilon}} v_\varepsilon^{- \frac{(2p-1)q}{2(p+1)}} |v_{\varepsilon x}|^q \right \}  \leq C(K) ||v_\varepsilon||_{L^\infty(B_{1/\varepsilon})} \Big( \int_{B_{1/\varepsilon}} u_\varepsilon^p  + 1 + \varepsilon^{\frac{p(q+2)+2}{q}} + \varepsilon^{\frac{q}{2}} \Big).
    $$
    Thus, as $\varepsilon < 1$, by defining
    $$ y_\varepsilon(t) := \int_{B_{1/\varepsilon}} u_\varepsilon^p(\cdot,t) + \int_{B_{1/\varepsilon}} v_\varepsilon^{- \frac{(2p-1)q}{2(p+1)}}(\cdot,t)  |v_{\varepsilon x}(\cdot,t)|^q + 3, \quad t \in [0, T_{\max, \varepsilon}),
    $$
    we obtain
    $$
    y'_\varepsilon(t) \leq C(K) ||v_\varepsilon||_{L^\infty(B_{1/\varepsilon})} \cdot  y_\varepsilon(t).
    $$
    Integrating the inequality yields for all $t \in (0, T_{\max, \varepsilon})$
    \begin{equation}\label{3.20}
    \begin{split}
    y_\varepsilon(t) \leq y_\varepsilon(0) \cdot e^{\displaystyle C(K) \int_0^t ||v_\varepsilon(\cdot,s)||_{L^\infty(B_{1/\varepsilon})}~ds}  \leq 
      y_\varepsilon(0) \cdot e^{\displaystyle C(K) ||v_0||_{L^\infty(\mathbb{R})} \cdot t}~,
    \end{split}
    \end{equation}
    where we used that by Lemma \ref{l2.2}, $0 \leq||v_\varepsilon(\cdot,t)||_{L^\infty(B_{1/\varepsilon})} \leq ||v_0||_{L^\infty(\mathbb{R})}$ for all $t \in (0, T_{\max, \varepsilon})$.
    \\\\
    Taking $C(p,K,t) : =  y_\varepsilon(0) \cdot e^{\displaystyle C(K) ||v_0||_{L^\infty(\mathbb{R})} \cdot t}$ finishes the proof.
\end{proof}
As a direct consequence of the $L^p$ bounds for $u_\varepsilon$, well-known properties of the Neumann heat semigroup can be applied to obtain the following result for $v_\varepsilon$. 
\begin{lemma}\label{l3.8}
    Let $K>0$, then for all $\varepsilon \in (0,1)$, $t \in (0,T_{\max, \varepsilon})$, there exists $C(K,t) >0$ independent of $\varepsilon$ such that if $u_0$ and $v_0$ satisfy \eqref{hip0}, \eqref{hip1} and \eqref{3.19}, then 
    $$
    ||v_{\varepsilon x}||_{L^\infty (B_{1/\varepsilon})} \leq C(K,t).
    $$
\end{lemma}
\begin{proof}
    The result is a consequence of standard semigroup theory based on the $L^p$ bound for $u_\varepsilon$ derived in Lemma \ref{l3.7} (which, although time dependent, is always finite) and the $W^{1,\infty}(\mathbb{R})$ bound for $v_0$. A detailed proof can be found in Lemma 2.2 in \cite{winkler2d}.
\end{proof}
As a second consequence of the $L^p$ estimates for $u_\varepsilon$, we can prove that indeed $T_{\max, \varepsilon} = \infty$ and therefore the regularized solutions do exist globally in time.
\begin{lemma}\label{l3.9}
    Let $K>0$ and assume that $u_0$ and $v_0$ are such that \eqref{hip0}, \eqref{hip1} and \eqref{3.19} hold. Then $T_{\max, \varepsilon} = \infty$ for all $\varepsilon \in (0,1)$.
\end{lemma}
\begin{proof}
    Again, the proof is standard and relies on the fact that if for any $\varepsilon \in (0,1)$, $T_{\max, \varepsilon}$ was finite, then for all $p \geq 2$, the quantity $ \sup_{t \in (0,T_{\max, \varepsilon})} ||u_{\varepsilon} (\cdot, t)||_{L^p(B_{1/\varepsilon})}$ would also be finite. In this case, we refer the reader to Lemma 4.1 in \cite{li-wink} for more details.
\end{proof}
Once the global existence of regularized solutions has been established on each ball $B_{1/\varepsilon}$ for all $\varepsilon \in (0,1)$, we derive further estimates uniform in $\varepsilon$ to pass to the limit in order to construct a global solution to the original problem in the whole space.

\section{Further temporal estimates: a bound for $(u_\varepsilon^{\frac{p+1}{2}} v_\varepsilon)_{\varepsilon \in (0,1)}$ in $L^2\big((0,T); W^{1,1}_\text{loc}(\mathbb{R})\big)$}\label{s4}
In this section, we obtain further time-dependent estimates with the aim of proving compactness properties of the sequence of regularized solutions, which will allow us to extract a converging subsequence by means of an Aubin-Lions type lemma. 
\\\\
First, given the local boundedness of $u_\varepsilon$ in $L^p(B_{1/\varepsilon})$, the nonlinear term $-u_\varepsilon v_\varepsilon$ in the second equation of the regularized system \eqref{reg} can be easily handled. This will result in the necessary features of $(v_\varepsilon)_{\varepsilon \in (0,1)}$ to grant the existence of a converging subsequence.
\\\\
With respect to $u_\varepsilon$, we consider the auxiliary sequence $\displaystyle (u_\varepsilon^{\frac{p+1}{2}} v_\varepsilon)_{\varepsilon\in(0,1)}$. In order to obtain boundedness in suitable spaces, we start by a first technical result which will be of key importance on the later analysis.
\begin{lemma}\label{l4.1}
    Let $K>0$ be such that \eqref{hip0}, \eqref{hip1} and \eqref{3.19} are satisfied, as well as
    \begin{equation}\label{4.1}
        \int_\mathbb{R} \frac{v_{0x}^2}{v_0} < K.
    \end{equation}
    Then, for all $T > 0$ there exists $C(T)>0$ such that
    $$\int_0^T \int_{B_{1/\varepsilon}} \frac{v_{\varepsilon x}^4}{v_\varepsilon^3} < C(T), \quad t \in (0,T), $$
    for all $\varepsilon \in (0,1)$.
\end{lemma}
\begin{proof}
By the positivity of $v_\varepsilon$, we can compute 
\begin{equation}\label{4.2}
        \frac{d}{dt} \int_{B_{1/\varepsilon}} \frac{v_{\varepsilon x}^2}{v_\varepsilon} = 2 \int_{B_{1/\varepsilon}} \frac{v_{\varepsilon x}}{v_\varepsilon} (v_{\varepsilon x})_t - \int_{B_{1/\varepsilon}} \frac{v_{\varepsilon x}^2}{v_\varepsilon^2} v_{\varepsilon t}, \quad \text{for all }t >0.
\end{equation}
 By standard parabolic theory, due to the regularity of $u_\varepsilon$ and $v_\varepsilon$ provided by Lemma \ref{l2.1}, we have that $v_{\varepsilon x} \in C^{2,1}(\bar{B}_{1/\varepsilon} \times (0, T_{\max, \varepsilon})$ and satisfies the differentiated version of the second equation in \eqref{reg}, this is
    $$ (v_{\varepsilon x})_t = v_{\varepsilon xxx} - u_{\varepsilon x} v_\varepsilon - u_\varepsilon v_{\varepsilon x}.$$
Thus, substituting in \eqref{4.2} and integrating by parts we obtain
\begin{equation}\label{4.3}
    \begin{split}
        & \frac{d}{dt} \int_{B_{1/\varepsilon}} \frac{v_{\varepsilon x}^2}{v_\varepsilon} = 2 \int_{B_{1/\varepsilon}} \frac{v_{\varepsilon x}}{v_\varepsilon} (v_{\varepsilon xxx} - u_{\varepsilon x} v_\varepsilon - u_\varepsilon v_{\varepsilon x}) - \int_{B_{1/\varepsilon}} \frac{v_{\varepsilon x}^2}{v_\varepsilon^2} (v_{\varepsilon xx} - u_\varepsilon v_\varepsilon) \\\\
        & = - 2 \int_{B_{1/\varepsilon}} \frac{v_{\varepsilon xx}^2}{v_\varepsilon} + 2\int_{B_{1/\varepsilon}} \frac{v_{\varepsilon x}^2}{v_\varepsilon^2} v_{\varepsilon xx} - 2 \int_{B_{1/\varepsilon}} u_{\varepsilon x} v_{\varepsilon x} - 2 \int_{B_{1/\varepsilon}} \frac{u_\varepsilon}{v_\varepsilon} v_{\varepsilon x}^2 - \int_{B_{1/\varepsilon}} \frac{v_{\varepsilon x}^2}{v_\varepsilon^2} v_{\varepsilon xx} + \int_{B_{1/\varepsilon}} \frac{u_\varepsilon}{v_\varepsilon} v_{\varepsilon x}^2
        \\\\
        & = - 2 \int_{B_{1/\varepsilon}} \frac{v_{\varepsilon xx}^2}{v_\varepsilon} + \int_{B_{1/\varepsilon}} \frac{v_{\varepsilon x}^2}{v_\varepsilon^2} v_{\varepsilon xx} - 2 \int_{B_{1/\varepsilon}} u_{\varepsilon x} v_{\varepsilon x} - \int_{B_{1/\varepsilon}} \frac{u_\varepsilon}{v_\varepsilon} v_{\varepsilon x}^2, \quad \text{for all }t >0.
    \end{split}
\end{equation}
To estimate the terms in the right hand side, for all $t>0$, firstly integrating by parts one obtains
\begin{equation}\label{4.4}
    \int_{B_{1/\varepsilon}} \frac{v_{\varepsilon x^2}}{v_\varepsilon^2} v_{\varepsilon xx} = \frac{2}{3} \int_{B_{1/\varepsilon}} \frac{v_{\varepsilon x}^4}{v_\varepsilon^3},
\end{equation}
and moreover, using Lemma 3.2 in \cite{winkler1d} we can prove that
\begin{equation}\label{4.5}
    - \int_{B_{1/\varepsilon}} \frac{v_{\varepsilon xx}^2}{v_\varepsilon} \leq -\frac{4}{9} \int_{B_{1/\varepsilon}} \frac{v_{\varepsilon x}^4}{v_\varepsilon^3}.
\end{equation}
Lastly, integrating first by parts and then using Young's inequality, we get
\begin{equation}\label{4.6}
    -2 \int_{B_{1/\varepsilon}} u_{\varepsilon x} v_{\varepsilon x} = 2 \int_{B_{1/\varepsilon}} u_\varepsilon v_{\varepsilon xx} \leq \frac{1}{4} \int_{B_{1/\varepsilon}} \frac{v_{\varepsilon xx}^2}{v_\varepsilon} + 4 \int_{B_{1/\varepsilon}} u_{\varepsilon}^2 v_\varepsilon.
\end{equation}
Substituting bounds \eqref{4.4}-\eqref{4.6} in \eqref{4.3} and using the positivity of $u_\varepsilon$ and $v_\varepsilon$ yields
\begin{equation}\label{4.7}
    \begin{split}
        \frac{d}{dt} \int_{B_{1/\varepsilon}} &\frac{v_{\varepsilon x}^2}{v_\varepsilon}  =   \left(-2 + \frac{1}{4} \right) \int_{B_{1/\varepsilon}} \frac{v_{\varepsilon xx}^2}{v_\varepsilon} + \int_{B_{1/\varepsilon}} \frac{v_{\varepsilon x}^2}{v_\varepsilon^2} v_{\varepsilon xx} + 4 \int_{B_{1/\varepsilon}} u_{\varepsilon}^2 v_\varepsilon - \int_{B_{1/\varepsilon}} \frac{u_\varepsilon}{v_\varepsilon} v_{\varepsilon x}^2
        \\\\
        & \leq - \frac{7}{4} \int_{B_{1/\varepsilon}} \frac{v_{\varepsilon xx}^2}{v_\varepsilon} + \frac{2}{3} \int_{B_{1/\varepsilon}} \frac{v_{\varepsilon x^4}}{v_\varepsilon^3} + 4 \int_{B_{1/\varepsilon}} u_{\varepsilon}^2 v_\varepsilon \leq \left(-\frac{7}{9} + \frac{2}{3} \right)  \int_{B_{1/\varepsilon}} \frac{v_{\varepsilon x}^4}{v_\varepsilon^3} + 4 \int_{B_{1/\varepsilon}} u_\varepsilon^2 v_\varepsilon,
    \end{split}
\end{equation}
for all $t >0$.
\\\\
In this way, for any given $T>0$, by the $L^\infty(B_{1/\varepsilon})$ bound for $v_\varepsilon$ from Lemma \ref{l2.1} and the $L^p$ bound for $u_\varepsilon$ provided by Lemma \ref{l3.7} there exists $C_1(T)>0$ such that
$$
4 \int_{B_{1/\varepsilon}} u_\varepsilon^2 v_\varepsilon \leq 4 ||v_0||_{L^\infty(\mathbb{R})} \int_{B_{1/\varepsilon}} u_\varepsilon^2 \leq C_1(T), \quad \text{for all } t \in (0,T),
$$
after which \eqref{4.7} can be rewritten as
\begin{equation} \label{4.8}
    \frac{d}{dt} \int_{B_{1/\varepsilon}} \frac{v_{\varepsilon x}^2}{v_\varepsilon} + \frac{1}{9}  \int_{B_{1/\varepsilon}} \frac{v_{\varepsilon x^4}}{v_\varepsilon^3} \leq C_1(T), \quad \text{for all } t \in (0,T).
\end{equation}
Integrating with respect to time provides
$$
\int_{B_{1/\varepsilon}} \frac{v_{\varepsilon x}^2}{v_\varepsilon} + \frac{1}{9} \int_0^T \int_{B_{1/\varepsilon}} \frac{v_{\varepsilon x}^4}{v_\varepsilon^3}  \leq \tilde{C}_1(T) + \int_{B_{1/\varepsilon}} \frac{v_{0x}^2}{v_0} \leq \tilde{C}_1(T) + K.
$$
by assumption \eqref{4.1}. Thus, the positivity of the first integral on the left hand side finishes the proof.
\end{proof}
Next, to derive the desired compactness properties of the sequence $(u_\varepsilon^{\frac{p+1}{2}} v_\varepsilon)_{\varepsilon\in(0,1)}$, we consider the following lemmas that will allow us to obtain local boundedness in $L^2\big((0,T); W^{1,1}_\text{loc}(\mathbb{R})\big)$. 
\\\\
In order to obtain estimates that remain uniform with respect to the domain size, we consider a general smooth, nonnegative, and compactly supported cutoff function as defined below which will allow us to localize the analysis.
\begin{definition}\label{ctf}
    Let $\varepsilon\in (0,1)$. We define a cutoff function $\phi^2 = \phi^2(x) \in C_c^\infty(B_{1/\varepsilon})$ such that
    \begin{itemize}
  \item $0\leq \phi^2 \leq 1$ in $B_{1/\varepsilon}$,
  \item $\phi^2 \equiv 1$ on a fixed compact subset $(-R, R)\subset\subset B_{1/\varepsilon}$,
  \item $\phi^2 \equiv 0$ in a neighborhood of $\partial B_{1/\varepsilon}$,
  \item $\phi^2$ connects smoothly the values $1$ in $(-R,R)$ and $0$ near $\partial B_{1/\varepsilon}$.
    \end{itemize}
\end{definition}
For such a class of functions, we prove the following results.
\begin{lemma}\label{l4.2}
    Let $p \geq 2$ and $K>0$ be such that \eqref{hip0}, \eqref{hip1} and \eqref{3.19} are satisfied. Then for any $\phi^2$ as in Definition \ref{ctf} and $T>0$ there exists $C(p,T)>0$ such that following inequality holds for all $\varepsilon \in (0,1)$
    \begin{equation*}
    \begin{split}
        \frac{d}{dt} \int_{B_{1/\varepsilon}} & u_\varepsilon^p \phi^2  + \frac{p (p-1)}{4} \int_{B_{1/\varepsilon}} u_\varepsilon^{p-1} v_\varepsilon u_{\varepsilon x}^2 \phi^2 + \frac{1}{2} \int_{B_{1/\varepsilon}} u_\varepsilon^2 v_\varepsilon \phi^2 \leq C(p,T) \int_{B_{1/\varepsilon}} u_\varepsilon^{p+1} v_\varepsilon \phi^2\\
        &  + \frac{8p}{p-1} \int_{B_{1/\varepsilon}} u_\varepsilon^{p+1} v_\varepsilon \phi_x^2  + C(p,T) \int_{B_{1/\varepsilon}} u_\varepsilon v_\varepsilon \phi^2, \quad \text{for all } t \in (0,T).
        \end{split}
\end{equation*}
\end{lemma}
\begin{proof}
    Using the first equation in \eqref{reg}, as $\phi^2 = \phi^2(x)$, and integrating by parts we obtain
 \begin{equation*}
     \begin{split}
         \frac{d}{dt} \int_{B_{1/\varepsilon}}  u_\varepsilon^p \phi^2 = \int_{B_{1/\varepsilon}} (u_\varepsilon^p)_t ~\phi^2 =  p \int_{B_{1/\varepsilon}} u_\varepsilon^{p-1} \Big((u_\varepsilon v_\varepsilon u_{\varepsilon x})_x - (u_\varepsilon^2 v_\varepsilon v_{\varepsilon x})_x + u_\varepsilon v_\varepsilon\Big) \phi^2 \\
          -p (p-1) \int_{B_{1/\varepsilon}} u_\varepsilon^{p-1} v_\varepsilon u_{\varepsilon x}^2 \phi^2 - 2p \int_{B_{1/\varepsilon}} u_\varepsilon^p v_\varepsilon u_{\varepsilon x} \phi \phi_x + p(p-1) \int_{B_{1/\varepsilon}} u_\varepsilon^p v_\varepsilon u_{\varepsilon x} v_{\varepsilon x} \phi^2 \\
          + 2p \int_{B_{1/\varepsilon}} u_\varepsilon^{p+1} v_\varepsilon v_{\varepsilon x} \phi \phi_x + p \int_{B_{1/\varepsilon}} u_\varepsilon^p v_\varepsilon \phi^2, \quad \text{for all } t \in (0,T).
     \end{split}
 \end{equation*}
    Thus, it follows directly that for all $t \in (0,T)$
    \begin{equation}\label{4.9}
        \begin{split}
             \frac{d}{dt} \int_{B_{1/\varepsilon}} & u_\varepsilon^p \phi^2  + p (p-1) \int_{B_{1/\varepsilon}} u_\varepsilon^{p-1} v_\varepsilon u_{\varepsilon x}^2 \phi^2 = - 2p \int_{B_{1/\varepsilon}} u_\varepsilon^p v_\varepsilon u_{\varepsilon x} \phi \phi_x + p(p-1) \int_{B_{1/\varepsilon}} u_\varepsilon^p v_\varepsilon u_{\varepsilon x} v_{\varepsilon x} \phi^2 \\
          &+ 2p \int_{B_{1/\varepsilon}} u_\varepsilon^{p+1} v_\varepsilon v_{\varepsilon x} \phi \phi_x + p \int_{B_{1/\varepsilon}} u_\varepsilon^p v_\varepsilon \phi^2.
        \end{split}
    \end{equation}
We estimate the terms appearing in the right hand side using Young's inequality as follows
\begin{equation}\label{4.10}
    - 2p \int_{B_{1/\varepsilon}} u_\varepsilon^p v_\varepsilon u_{\varepsilon x} \phi \phi_x \leq \frac{p(p-1)}{4} \int_{B_{1/\varepsilon}} u_\varepsilon^{p-1} v_\varepsilon u_{\varepsilon x}^2 \phi^2 + \frac{4p}{p-1} \int_{B_{1/\varepsilon}} u_\varepsilon^{p+1} v_\varepsilon \phi_x^2,  \quad \text{for all } t \in (0,T),
\end{equation}
as well as
\begin{equation}\label{4.11}
\begin{split}
    p(p-1) \int_{B_{1/\varepsilon}} u_\varepsilon^p v_\varepsilon u_{\varepsilon x} v_{\varepsilon x} \phi^2 \leq \frac{p(p-1)}{2} \int_{B_{1/\varepsilon}} u_\varepsilon^{p-1} v_\varepsilon u_{\varepsilon x}^2 \phi^2 + \frac{p(p-1)}{2} \int_{B_{1/\varepsilon}} u_\varepsilon^{p+1} v_\varepsilon v_{\varepsilon x}^2 \phi^2 \\\\
    \leq \frac{p(p-1)}{2} \int_{B_{1/\varepsilon}} u_\varepsilon^{p-1} v_\varepsilon u_{\varepsilon x}^2 \phi^2 + \frac{p(p-1)}{2} c_1^2(T)  \int_{B_{1/\varepsilon}} u_\varepsilon^{p+1} v_\varepsilon \phi^2  \quad \text{for all } t \in (0,T),
\end{split}
\end{equation}
where we used that by Lemma \ref{l3.8} there exists $c_1(T) >0$ such that $||v_{\varepsilon x}||_{L^\infty (B_{1/\varepsilon})} \leq c_1(T)$. Again by the same argument we have
\begin{equation}\label{4.12}
\begin{split}
  2p \int_{B_{1/\varepsilon}} u_\varepsilon^{p+1} v_\varepsilon v_{\varepsilon x} \phi \phi_x \leq \frac{p(p-1)}{4} \int_{B_{1/\varepsilon}} u_\varepsilon^{p+1} v_\varepsilon v_{\varepsilon x}^2 \phi^2 + \frac{4p}{p-1} \int_{B_{1/\varepsilon}} u_\varepsilon^{p+1} v_\varepsilon \phi_x^2  \\\\
  \leq \frac{p(p-1)}{4} c_1^2(T) \int_{B_{1/\varepsilon}} u_\varepsilon^{p+1} v_\varepsilon  \phi^2 + \frac{4p}{p-1} \int_{B_{1/\varepsilon}} u_\varepsilon^{p+1} v_\varepsilon \phi_x^2 
   \quad \text{for all } t \in (0,T).
\end{split}
\end{equation}
Substituting \eqref{4.10}-\eqref{4.12} in \eqref{4.9} one obtains for all $t \in (0,T)$
\begin{equation}\label{4.13}
     \begin{split}
        \frac{d}{dt} \int_{B_{1/\varepsilon}} & u_\varepsilon^p \phi^2  + p (p-1) \int_{B_{1/\varepsilon}} u_\varepsilon^{p-1} v_\varepsilon u_{\varepsilon x}^2 \phi^2 \leq \frac{3p(p-1)}{4} \int_{B_{1/\varepsilon}} u_\varepsilon^{p-1} v_\varepsilon u_{\varepsilon x}^2 \phi^2  \\\\
        & + \frac{3p(p-1)}{4}c_1^2(T) \int_{B_{1/\varepsilon}}  u_\varepsilon^{p+1} v_\varepsilon  \phi^2 + \frac{8p}{p-1} \int_{B_{1/\varepsilon}} u_\varepsilon^{p+1} v_\varepsilon \phi_x^2 + p \int_{B_{1/\varepsilon}} u_\varepsilon^p v_\varepsilon \phi^2.
    \end{split}
\end{equation}
For the last term, which we are yet to estimate, Young's inequality, this time with exponents $\displaystyle \frac{p}{p-1}$ and $p$---both greater than 1 as $p \geq 2$---yields
\begin{equation}\label{4.14}
    \begin{split}
        p \int_{B_{1/\varepsilon}} u_\varepsilon^p v_\varepsilon \phi^2 &= \int_{B_{1/\varepsilon}} \left( \frac{3p}{4} c_1^2(T) \cdot \frac{p}{p-1} u_\varepsilon^{p+1} v_\varepsilon \right)^{\frac{p-1}{p}} \cdot \left( p \left( \frac{4(p-1)}{3 p^2 c_1^2(T)} \right)^{\frac{p-1}{p}} u_\varepsilon^{1/p} v_\varepsilon^{1/p} \right) \phi^2 \\\\
        & \leq \frac{3p}{4} c_1^2(T) \int_{B_{1/\varepsilon}}  u_\varepsilon^{p+1} v_\varepsilon \phi^2 + \frac{1}{p} \left \{ p \left( \frac{4(p-1)}{3 p^2 c_1^2(T)} \right)^{\frac{p-1}{p}}\right\}^p \int_{B_{1/\varepsilon}} u_\varepsilon v_\varepsilon \phi^2, \quad \text{for all } t \in (0,T).
    \end{split}
\end{equation}
In this way, the first summand on the right hand side here cancels the same term arising from the negative contribution in $\displaystyle \frac{3p(p-1)}{4}c_1^2(T)$ in \eqref{4.13}.
\\\\
Lastly, following the same strategy, we obtain for all $ t \in (0,T)$
\begin{equation}\label{4.15}
    \begin{split}
        \frac{1}{2} \int_{B_{1/\varepsilon}} u_\varepsilon^2 v_\varepsilon \phi^2 & = \int_{B_{1/\varepsilon}} \left( \frac{p^2}{4} c_1^2(T) \cdot p~ u_\varepsilon^{p+1} v_\varepsilon \right)^{1/p} \cdot \left( \frac{1}{2} \cdot \left(\frac{4}{p^3 c_1^2(T)} \right)^{1/p} (u_\varepsilon v_\varepsilon)^{\frac{p-1}{p}} \right) \phi^2 \\\\
        & \leq \frac{p^2}{4} c_1^2(T) \int_{B_{1/\varepsilon}} u_\varepsilon^{p+1} v_\varepsilon \phi^2  + \frac{p-1}{p} \left\{\frac{1}{2} \cdot \left(\frac{4}{p^3 c_1^2(T)} \right)^{1/p} \right\}^{\frac{p}{p-1}} \int_{B_{1/\varepsilon}} u_\varepsilon v_\varepsilon \phi^2.
    \end{split}
\end{equation}
Therefore, combining \eqref{4.14} and \eqref{4.15} into \eqref{4.13} yields
\begin{equation*}
    \begin{split}
        \frac{d}{dt} \int_{B_{1/\varepsilon}} & u_\varepsilon^p \phi^2  + \frac{p (p-1)}{4} \int_{B_{1/\varepsilon}} u_\varepsilon^{p-1} v_\varepsilon u_{\varepsilon x}^2 \phi^2 + \frac{1}{2} \int_{B_{1/\varepsilon}} u_\varepsilon^2 v_\varepsilon \phi^2 \leq p^2 c_1^2(T) \int_{B_{1/\varepsilon}} u_\varepsilon^{p+1} v_\varepsilon \phi^2\\
        &  + \frac{8p}{p-1} \int_{B_{1/\varepsilon}} u_\varepsilon^{p+1} v_\varepsilon \phi_x^2  + \left[ \left( \frac{4(p-1)}{3 p c_1^2(T)} \right)^{p-1} + \frac{p-1}{p} \cdot \frac{1}{2^{\frac{p}{p-1}}} \left(\frac{4}{p^3 c_1^2(T)} \right)^{\frac{1}{p-1}} \right] \int_{B_{1/\varepsilon}} u_\varepsilon v_\varepsilon \phi^2,
        \end{split}
\end{equation*}
which gives the result considering $C(p,T)$ as
$$
C(p,T) := \max \left\{p^2 c_1^2(T) , \displaystyle \left[ \left( \frac{4(p-1)}{3 p c_1^2(T)} \right)^{p-1} + \frac{p-1}{p} \cdot \frac{1}{2^{\frac{p}{p-1}}} \left(\frac{4}{p^3 c_1^2(T)} \right)^{\frac{1}{p-1}} \right] \right \}>0.
$$
\end{proof}
Next, we study the time derivative of $\displaystyle \int_{B_{1/\varepsilon}} \frac{|v_{\varepsilon x}|^2}{v_\varepsilon} \phi^2$ in the following Lemma.
\begin{lemma}\label{l4.3}
Let $K>0$ such that \eqref{hip0}, \eqref{hip1} and \eqref{3.19} are satisfied. Then for any $\phi^2$ as in Definition \ref{ctf} we have 
    \begin{equation*}
    \begin{split}
        \frac{d}{dt} \int_{B_{1/\varepsilon}} &\frac{v_{\varepsilon x}^2}{v_\varepsilon} \phi^2  + \int_{B_{1/\varepsilon}}\frac{u_\varepsilon}{v_{\varepsilon}}v_{\varepsilon x}^2 \phi^2  \leq \frac{1}{2} \int_{B_{1/\varepsilon}} u_\varepsilon^2 v_\varepsilon \phi^2 + 2 \int_{B_{1/\varepsilon}} u_\varepsilon v_{\varepsilon x} (\phi^2)_x \\\\
            &+4 \int_{B_{1/\varepsilon}} \frac{v_{\varepsilon x}^2 v_{\varepsilon xx}}{v_\varepsilon^2} \phi^2  + \int_{B_{1/\varepsilon}} \frac{v_{\varepsilon x}^2}{v_\varepsilon} (\phi^2)_{xx}, \quad \text{for all } t >0,
        \end{split}
    \end{equation*}
    for any choice of $\varepsilon \in (0,1)$.
\end{lemma}
\begin{proof}
    Computing the time derivative, we have
    \begin{equation*}
        \begin{split}
            \frac{d}{dt}&\int_{B_{1/\varepsilon}} \frac{v_{\varepsilon x}^2}{v_\varepsilon} \phi^2 = \int_{B_{1/\varepsilon}} \frac{2v_{\varepsilon x} v_{\varepsilon xt} v_\varepsilon - v_{\varepsilon x}^2 v_{\varepsilon t}}{v_\varepsilon^2}\phi^2 = 2 \int_{B_{1/\varepsilon}} \frac{ v_{\varepsilon x} (v_{\varepsilon xx} - u_\varepsilon v_\varepsilon)_x}{v_\varepsilon}\phi^2 \\\\
            &- \int_{B_{1/\varepsilon}} \frac{v_{\varepsilon x}^2(v_{\varepsilon xx}-u_\varepsilon v_\varepsilon)}{v_\varepsilon^2} \phi^2 = 2\int_{B_{1/\varepsilon}} \frac{ v_{\varepsilon x} v_{\varepsilon xxx}}{v_\varepsilon} \phi^2 - 2\int_{B_{1/\varepsilon}} u_{\varepsilon x} v_{\varepsilon x} \phi^2 - 2 \int_{B_{1/\varepsilon}}\frac{u_\varepsilon}{v_{\varepsilon}}v_{\varepsilon x}^2 \phi^2 \\\\
            & - \int_{B_{1/\varepsilon}} \frac{v_{\varepsilon x}^2 v_{\varepsilon xx}}{v_\varepsilon^2} \phi^2+ \int_{B_{1/\varepsilon}} \frac{u_\varepsilon}{v_{\varepsilon}}v_{\varepsilon x}^2 \phi^2, \quad \text{for all } t>0.
            \end{split}
    \end{equation*}
    Rewriting the expression, we have for all $t >0$ 
    \begin{equation}\label{4.16}
        \begin{split}
            \frac{d}{dt} \int_{B_{1/\varepsilon}} \frac{v_{\varepsilon x}^2}{v_\varepsilon} \phi^2  + \int_{B_{1/\varepsilon}}\frac{u_\varepsilon}{v_{\varepsilon}}v_{\varepsilon x}^2 \phi^2  = - 2 \int_{B_{1/\varepsilon}} u_{\varepsilon x} v_{\varepsilon x} \phi^2 + 
            2 \int_{B_{1/\varepsilon}} \frac{v_{\varepsilon x} v_{\varepsilon xxx}}{v_\varepsilon} \phi^2 - \int_{B_{1/\varepsilon}} \frac{v_{\varepsilon x}^2 v_{\varepsilon xx}}{v_\varepsilon^2} \phi^2.
        \end{split}
    \end{equation}
    Next, first integrating by parts and then using Young's inequality we have
    \begin{equation}\label{4.17}
        \begin{split}
            - 2 \int_{B_{1/\varepsilon}} & u_{\varepsilon x} v_{\varepsilon x} \phi^2 = 2 \int_{B_{1/\varepsilon}} u_\varepsilon v_{\varepsilon xx} \phi^2 + 2 \int_{B_{1/\varepsilon}} u_\varepsilon v_{\varepsilon x} (\phi^2)_x \\
            & \leq 2 \int_{B_{1/\varepsilon}} \frac{v_{\varepsilon xx}^2}{v_\varepsilon} \phi^2 + \frac{1}{2} \int_{B_{1/\varepsilon}} u_\varepsilon^2 v_\varepsilon \phi^2 + 2 \int_{B_{1/\varepsilon}} u_\varepsilon v_{\varepsilon x} (\phi^2)_x ,
        \end{split}
    \end{equation}
    Now we deal with the term involving $v_{\varepsilon xxx}$. To do so, we make use of the following identity for any $f \in C^3(\mathbb{R})$
    $$
    \frac{d^2}{dx^2}\Big((f')^2\Big) = 2 (f'')^2 + 2 f' f''',
    $$
    which in particular for all $t >0$ implies
    $$
    2 v_{\varepsilon x} v_{\varepsilon xxx} = (v_{\varepsilon x}^2)_{xx} - 2 v_{\varepsilon xx}^2.
    $$
    Hence, after integrating by parts we obtain for all $t>0$
    \begin{equation}\label{4.17a}
      \begin{split}
          2 \int_{B_{1/\varepsilon}}& \frac{v_{\varepsilon x} v_{\varepsilon xxx}}{v_\varepsilon} \phi^2 = \int_{B_{1/\varepsilon}} \frac{ (v_{\varepsilon x}^2)_{xx}}{v_\varepsilon} \phi^2 - 2 \int_{B_{1/\varepsilon}} \frac{v_{\varepsilon xx}^2}{v_\varepsilon}\phi^2 = \int_{B_{1/\varepsilon}} v_{\varepsilon x}^2 \left( \frac{\phi^2}{v_\varepsilon}\right)_{xx} \hspace{-0.2 cm} - 2 \int_{B_{1/\varepsilon}} \frac{v_{\varepsilon xx}^2}{v_\varepsilon}\phi^2 \\\\
          & = \int_{B_{1/\varepsilon}} v_{\varepsilon x}^2 \left(\frac{(\phi^2)_x}{v_\varepsilon} - \frac{v_{\varepsilon x}}{v_\varepsilon^2}\phi^2 \right)_x - 2 \int_{B_{1/\varepsilon}} \frac{v_{\varepsilon xx}^2}{v_\varepsilon}\phi^2 = \int_{B_{1/\varepsilon}} \frac{v_{\varepsilon x}^2}{v_\varepsilon} (\phi^2)_{xx} - \int_{B_{1/\varepsilon}} \frac{v_{\varepsilon x}^3}{v_\varepsilon ^2} (\phi^2)_x \\\\
          &  - \int_{B_{1/\varepsilon}} \frac{v_{\varepsilon x}^2 v_{\varepsilon xx}}{v_\varepsilon^2} \phi^2 + 2 \int_{B_{1/\varepsilon}} \frac{v_{\varepsilon x}^4}{v_\varepsilon^3}\phi^2 - \int_{B_{1/\varepsilon}} \frac{v_{\varepsilon x}^3}{v_\varepsilon^2} (\phi^2)_x - 2 \int_{B_{1/\varepsilon}} \frac{v_{\varepsilon xx}^2}{v_\varepsilon}\phi^2 \\\\
          & = - \int_{B_{1/\varepsilon}} \frac{v_{\varepsilon x}^2 v_{\varepsilon xx}}{v_\varepsilon^2} \phi^2  + 2 \int_{B_{1/\varepsilon}} \frac{v_{\varepsilon x}^4}{v_\varepsilon^3}\phi^2 - 2 \int_{B_{1/\varepsilon}} \frac{v_{\varepsilon x}^3}{v_\varepsilon^2} (\phi^2)_x + \int_{B_{1/\varepsilon}} \frac{v_{\varepsilon x}^2}{v_\varepsilon} (\phi^2)_{xx} - 2 \int_{B_{1/\varepsilon}} \frac{v_{\varepsilon xx}^2}{v_\varepsilon}\phi^2.
      \end{split}  
   \end{equation}
   Notice that the last term precisely appears on the right hand side of \eqref{4.17} with opposite sign. One last integration by parts reveals
    $$
    -2 \int_{B_{1/\varepsilon}} \frac{v_{\varepsilon x}^3}{v_\varepsilon^2} (\phi^2)_x = 6 \int_{B_{1/\varepsilon}} \frac{v_{\varepsilon x}^2 v_{\varepsilon xx}}{v_\varepsilon^2} \phi^2 - 4 \int_{B_{1/\varepsilon}} \frac{v_{\varepsilon x}^4}{v_\varepsilon^3} \phi^2.
    $$
   Now, substituting this \eqref{4.17a} yields for all $t>0$
   \begin{equation}\label{4.18}
          2 \int_{B_{1/\varepsilon}} \frac{v_{\varepsilon x} v_{\varepsilon xxx}}{v_\varepsilon} \phi^2 = 5 \int_{B_{1/\varepsilon}} \frac{v_{\varepsilon x}^2 v_{\varepsilon xx}}{v_\varepsilon^2} \phi^2 - 2 \int_{B_{1/\varepsilon}} \frac{v_{\varepsilon x}^4}{v_\varepsilon^3} \phi^2 + \int_{B_{1/\varepsilon}} \frac{v_{\varepsilon x}^2}{v_\varepsilon} (\phi^2)_{xx} - 2 \int_{B_{1/\varepsilon}} \frac{v_{\varepsilon xx}^2}{v_\varepsilon}\phi^2,
   \end{equation}
   where the last term can be dropped due to its non-positivity. Thus, by \eqref{4.17} and \eqref{4.18}, \eqref{4.16} can be rewritten as
   \begin{equation*}
        \begin{split}
            \frac{d}{dt} \int_{B_{1/\varepsilon}} &\frac{v_{\varepsilon x}^2}{v_\varepsilon} \phi^2  + \int_{B_{1/\varepsilon}}\frac{u_\varepsilon}{v_{\varepsilon}}v_{\varepsilon x}^2 \phi^2  \leq \frac{1}{2} \int_{B_{1/\varepsilon}} u_\varepsilon^2 v_\varepsilon \phi^2 + 2 \int_{B_{1/\varepsilon}} u_\varepsilon v_{\varepsilon x} (\phi^2)_x \\\\
            &+4 \int_{B_{1/\varepsilon}} \frac{v_{\varepsilon x}^2 v_{\varepsilon xx}}{v_\varepsilon^2} \phi^2 + \int_{B_{1/\varepsilon}} \frac{v_{\varepsilon x}^2}{v_\varepsilon} (\phi^2)_{xx} , \quad \text{for all } t >0,
        \end{split}
   \end{equation*}    
    which finishes the proof.
\end{proof}
Next, we can readily join both previous lemmas to obtain integrability properties that will later grant adequate compactness of $(u_\varepsilon^{\frac{p+1}{2}} v_\varepsilon)_{\varepsilon\in(0,1)}$.
\begin{lemma}\label{l4.4}
    Let $p \geq 2$ and $K>0$ be such that \eqref{hip0}, \eqref{hip1}, \eqref{3.19} and \eqref{4.1} are satisfied. Then for any $\phi^2$ as in Definition \ref{ctf} and $T >0$, there exists $C(p,T) >0$ that verifies
    $$
    \int_0^T \int_{B_{1/\varepsilon}} u^{p-1} v_\varepsilon u_{\varepsilon x}^2 \phi^2 + \int_{B_{1/\varepsilon}} \frac{v_{\varepsilon x}(\cdot,T)^2}{v_\varepsilon (\cdot, T)} \phi^2 + \int_0^T \int_{B_{1/\varepsilon}} \frac{u_\varepsilon}{v_\varepsilon} v_{\varepsilon x}^2 \phi^2 \leq C(p,T),
    $$
    for all $\varepsilon \in (0,1)$.
\end{lemma}
\begin{proof}
    We begin by combining the results of Lemmas \ref{l4.2} and \ref{l4.3}. In this way, given any fixed $T>0$, for all $t \in (0, T)$ by Lemma \ref{l4.2} there exists $\hat{c}(p,T) >0$ such that
     \begin{equation*}
    \begin{split}
        \frac{d}{dt} \int_{B_{1/\varepsilon}} & u_\varepsilon^p \phi^2  + \frac{p (p-1)}{4} \int_{B_{1/\varepsilon}} u_\varepsilon^{p-1} v_\varepsilon u_{\varepsilon x}^2 \phi^2 + \frac{1}{2} \int_{B_{1/\varepsilon}} u_\varepsilon^2 v_\varepsilon \phi^2 \leq \hat{c}(p,T) \int_{B_{1/\varepsilon}} u_\varepsilon^{p+1} v_\varepsilon \phi^2\\
        &  + \frac{8p}{p-1} \int_{B_{1/\varepsilon}} u_\varepsilon^{p+1} v_\varepsilon \phi_x^2  + \hat{c}(p,T) \int_{B_{1/\varepsilon}} u_\varepsilon v_\varepsilon \phi^2, \quad \text{for all } t \in (0,T),
        \end{split}
\end{equation*}
    Similarly, Lemma \ref{l4.3} ensures that
    \begin{equation*}
    \begin{split}
        \frac{d}{dt} \int_{B_{1/\varepsilon}} &\frac{v_{\varepsilon x}^2}{v_\varepsilon} \phi^2  + \int_{B_{1/\varepsilon}}\frac{u_\varepsilon}{v_{\varepsilon}}v_{\varepsilon x}^2 \phi^2  \leq \frac{1}{2} \int_{B_{1/\varepsilon}} u_\varepsilon^2 v_\varepsilon \phi^2 + 2 \int_{B_{1/\varepsilon}} u_\varepsilon v_{\varepsilon x} (\phi^2)_x \\\\
            &+4 \int_{B_{1/\varepsilon}} \frac{v_{\varepsilon x}^2 v_{\varepsilon xx}}{v_\varepsilon^2} \phi^2  + \int_{B_{1/\varepsilon}} \frac{v_{\varepsilon x}^2}{v_\varepsilon} (\phi^2)_{xx} , \quad \text{for all } t \in (0,T).
        \end{split}
    \end{equation*}
    Adding both expressions we obtain that for all $t \in (0,T)$
    \begin{equation} \label{4.19}
        \begin{split}
          \frac{d}{dt} & \left \{\int_{B_{1/\varepsilon}} u_\varepsilon^p \phi^2 +  \int_{B_{1/\varepsilon}} \frac{v_{\varepsilon x}^2}{v_\varepsilon} \phi^2\right\} + \frac{p (p-1)}{4} \int_{B_{1/\varepsilon}} u_\varepsilon^{p-1} v_\varepsilon u_{\varepsilon x}^2 \phi^2 + \int_{B_{1/\varepsilon}}\frac{u_\varepsilon}{v_{\varepsilon}}v_{\varepsilon x}^2 \phi^2  \\\\
          & \leq \hat{c}(p,T) \int_{B_{1/\varepsilon}} u_\varepsilon^{p+1} v_\varepsilon \phi^2  + \frac{8p}{p-1} \int_{B_{1/\varepsilon}} u_\varepsilon^{p+1} v_\varepsilon \phi_x^2  + \hat{c}(p,T) \int_{B_{1/\varepsilon}} u_\varepsilon v_\varepsilon \phi^2 \\\\
          &+ 2 \int_{B_{1/\varepsilon}} u_\varepsilon v_{\varepsilon x} (\phi^2)_x 
            +4 \int_{B_{1/\varepsilon}} \frac{v_{\varepsilon x}^2 v_{\varepsilon xx}}{v_\varepsilon^2} \phi^2 + \int_{B_{1/\varepsilon}} \frac{v_{\varepsilon x}^2}{v_\varepsilon} (\phi^2)_{xx},
        \end{split}
    \end{equation}
    so a time integration provides
    \begin{equation}\label{4.20}
        \begin{split}
            \int_{B_{1/\varepsilon}} & u_\varepsilon^p(\cdot,T)~ \phi^2 +  \int_{B_{1/\varepsilon}} \frac{v_{\varepsilon x}(\cdot,T)^2}{v_\varepsilon(\cdot,T)} \phi^2 + \frac{p (p-1)}{4} \int_0^T \int_{B_{1/\varepsilon}} u_\varepsilon^{p-1} v_\varepsilon u_{\varepsilon x}^2 \phi^2 + \int_0^T \int_{B_{1/\varepsilon}}\frac{u_\varepsilon}{v_{\varepsilon}}v_{\varepsilon x}^2 \phi^2 \\\\\
            & \leq \int_{B_{1/\varepsilon}} (u_0 + \varepsilon \zeta)\phi^2 + \int_{B_{1/\varepsilon}} \frac{(v_0)_x^2}{v_0} \phi^2 + \hat{c}(p,T) \int_0^T \int_{B_{1/\varepsilon}} u_\varepsilon^{p+1} v_\varepsilon \phi^2  + \frac{8p}{p-1} \int_0^T  \int_{B_{1/\varepsilon}} u_\varepsilon^{p+1} v_\varepsilon \phi_x^2 \\\\\
            & +  \hat{c}(p,T)\int_0^T \int_{B_{1/\varepsilon}} u_\varepsilon v_\varepsilon \phi^2 + 2\int_0^T\int_{B_{1/\varepsilon}} u_\varepsilon v_{\varepsilon x} (\phi^2)_x 
            +4 \int_0^T \int_{B_{1/\varepsilon}} \frac{v_{\varepsilon x}^2 v_{\varepsilon xx}}{v_\varepsilon^2} \phi^2 \\\\
            & + \int_0^T \int_{B_{1/\varepsilon}} \frac{v_{\varepsilon x}^2}{v_\varepsilon} (\phi^2)_{xx} =: \int_{B_{1/\varepsilon}} (u_0 + \varepsilon \zeta)\phi^2 + \int_{B_{1/\varepsilon}} \frac{v_{0x}^2}{v_0} \phi^2 + \sum_{i=1}^6 I_i.
        \end{split}
    \end{equation}
    The first two terms on the right hand side are bounded by \eqref{hip1} and \eqref{4.1}, so to conclude the proof, we bound the other six terms involved.
    \\\\
    Firstly, due to Lemma \ref{l2.2} and the $L^p$ for arbitrary $p$ bounds proved in Lemma \ref{l3.7} we have
    \begin{equation*}
    \begin{split}
    I_1 = & \hat{c}(p,T) \int_0^T\int_{B_{1/\varepsilon}} u_\varepsilon^{p+1} v_\varepsilon \phi^2 \leq ||v_0||_{L^\infty(B_{1/\varepsilon})} \hat{c}(p,T) \int_0^T \int_{B_{1/\varepsilon}} u_\varepsilon^{p+1} \\
    & \leq ||v_0||_{L^\infty(\mathbb{R})} \hat{c}(p,T) \int_0^T c(p+1,T) ~dt := C_1(p,T) .
    \end{split}  
    \end{equation*}
    Next, with similar bounds and making use of Young's inequality
      \begin{equation*}
    \begin{split}
    I_2 = &\frac{8p}{p-1} \int_0^T  \int_{B_{1/\varepsilon}} u_\varepsilon^{p+1} v_\varepsilon \phi_x^2 \leq \frac{8p}{p-1} ||v_0||_{L^\infty(B_{1/\varepsilon})} \int_0^T  \int_{B_{1/\varepsilon}} u^{p+1} \phi_x^2 \\
    & \leq \frac{4p}{p-1} ||v_0||_{L^\infty(\mathbb{R})} \int_0^T  \int_{B_{1/\varepsilon}} \left(u_\varepsilon^{2(p+1)} + \phi_x^4 \right) \\
    & \leq \frac{4p}{p-1} ||v_0||_{L^\infty(\mathbb{R})} \int_0^T  \left ( c(2(p+1),T) +  \int_{B_{1/\varepsilon}} \phi_x^4 \right) ~dt := C_2(p,T).
    \end{split}
    \end{equation*}
    For $I_3$, the time integrability property of Lemma \ref{l2.2} directly provides
    $$
    I_3 = \hat{c}(p,T) \int_0^T \int_{B_{1/\varepsilon}} u_\varepsilon v_\varepsilon \phi^2 \leq \hat{c}(p,T) \int_0^T \int_{B_{1/\varepsilon}} u_\varepsilon v_\varepsilon \leq \hat{c}(p,T) \int_{\mathbb{R}} v_0 := C_3(p,T).
    $$
    With respect to $I_4$, considering the bound for $v_{\varepsilon x}$ in Lemma \ref{l3.8} there exists $\bar{c}(T)>0$ such that $\displaystyle ||v_{\varepsilon x}||_{L^\infty(B_{1/\varepsilon})} \leq \bar{c}(T)$. Hence, again using Young's inequality and the $L^p$ bounds provided by Lemma \ref{l3.7}
    \begin{equation*}
    \begin{split}
    I_4 &= 2\int_0^T\int_{B_{1/\varepsilon}} u_\varepsilon v_{\varepsilon x} (\phi^2)_x \leq 2 \bar{c}(T) \int_0^T \int_{B_{1/\varepsilon}} u_\varepsilon (\phi^2)_x \leq \bar{c}(T) \int_0^T \left( \int_{B_{1/\varepsilon}} u_\varepsilon^2 + \int_{B_{1/\varepsilon}}(\phi^2)_x \right) \\
    &\leq \bar{c}(T) \int_0^T \left( c(T) + \int_{B_{1/\varepsilon}}(\phi^2)_x \right)~dt  := C_4(T).    
    \end{split}
    \end{equation*}
    The analysis of $I_5$ is slightly more involved, and requires the key estimate provided by Lemma \ref{l4.1}. First, integrating by parts we have
    \begin{equation*}
        \begin{split}
           I(t):=&\int_{B_{1/\varepsilon}} \frac{v_{\varepsilon x}^2 v_{\varepsilon xx}}{v_\varepsilon^2} \phi^2 = - \int_{B_{1/\varepsilon}} v_{\varepsilon x} \left(\frac{v_{\varepsilon x}^2}{v_\varepsilon^2} \phi^2 \right)_x \\
           &= - 2 \int_{B_{1/\varepsilon}} \frac{v_{\varepsilon x}^2 v_{\varepsilon xx}}{v_\varepsilon^2} \phi^2 + 2 \int_{B_{1/\varepsilon}} \frac{v_{\varepsilon x}^4}{v_\varepsilon^3} \phi^2 - \int_{B_{1/\varepsilon}} \frac{v_{\varepsilon x}^3}{v_\varepsilon^2} (\phi^2)_x\\
           &= -2 I(t) + 2 \int_{B_{1/\varepsilon}} \frac{v_{\varepsilon x}^4}{v_\varepsilon^3} \phi^2 - \int_{B_{1/\varepsilon}} \frac{v_{\varepsilon x}^3}{v_\varepsilon^2} (\phi^2)_x
        \end{split}
    \end{equation*}
    Therefore we have 
    $$I_5 = \int_0^T I(t) ~dt = \frac{2}{3}\int_0^T  \int_{B_{1/\varepsilon}} \frac{v_{\varepsilon x}^4}{v_\varepsilon^3} \phi^2 - \frac{1}{3}\int_{B_{1/\varepsilon}}  \int_0^T\frac{v_{\varepsilon x}^3}{v_\varepsilon^2} (\phi^2)_x.$$
    As a direct consequence of Lemma \ref{l4.1}, the first term is indeed bounded, whereas for the second one, using Young's inequality with exponents $4/3$ and $4$ we obtain
    \begin{equation*}
        \begin{split}
            - \int_{B_{1/\varepsilon}} & \frac{v_{\varepsilon x}^3}{v_\varepsilon^2} (\phi^2)_x  \leq \int_{B_{1/\varepsilon}} \frac{|v_{\varepsilon x}|^3}{v_\varepsilon^2} |(\phi^2)_x| =  \int_{B_{1/\varepsilon}} \frac{|v_{\varepsilon x}|^3}{ v_\varepsilon^{9/4}} \cdot v_\varepsilon^{1/4} |(\phi^2)_x| \\
            & \leq \frac{3}{4} \int_{B_{1/\varepsilon}} \frac{|v_{\varepsilon x}|^4}{v_\varepsilon^3} + \frac{1}{4}\int_{B_{1/\varepsilon}} v_\varepsilon |(\phi^2)_x|^4 \leq \frac{3}{4} \int_{B_{1/\varepsilon}} \frac{v_{\varepsilon x}^4}{v_\varepsilon^3} + \frac{||v_0||_{L^\infty(\mathbb{R})}}{4} \int_{B_{1/\varepsilon}} |(\phi^2)_x|^4.
        \end{split}
    \end{equation*}
    The time integral of the first term is again bounded by Lemma \ref{l4.1}, while the second term is a constant independent of $\varepsilon$ that can be integrated, granting the existence of $C_5(T)$ such that $\displaystyle I_5 \leq C_5(T)$.
    \\\\
    Lastly, for $I_6$ a similar argument can be applied. Using Young's inequality combined with the $L^\infty$ bound for $v_\varepsilon$, one has
    \begin{equation*}
        \begin{split}
            I_6 &= \int_0^T\int_{B_{1/\varepsilon}} \frac{v_{\varepsilon x}^2}{v_\varepsilon} (\phi^2)_{xx} \leq \int_0^T\int_{B_{1/\varepsilon}} \frac{v_{\varepsilon x}^4}{v_\varepsilon^2} + \frac{1}{4} \int_0^T\int_{B_{1/\varepsilon}} \Big[(\phi^2)_{xx}\Big]^2 \\
            & \leq ||v_0||_{L^\infty(\mathbb{R})} \int_0^T\int_{B_{1/\varepsilon}} \frac{v_{\varepsilon x}^4}{v_\varepsilon^3}  + \frac{1}{4} \int_0^T\int_{B_{1/\varepsilon}} \Big[(\phi^2)_{xx}\Big]^2 \leq C_6(T),
        \end{split}
    \end{equation*}
    for a certain $C_6(T)>0$, thanks again to Lemma \ref{l4.1} that allows us to bound $\displaystyle \int_0^T\int_{B_{1/\varepsilon}} \frac{v_{\varepsilon x}^4}{v_\varepsilon^3}$.
    \\\\
    Thus, combining the bounds ensures the existence of $C(p,T)>0$ such that in particular
    $$
    \int_0^T \int_{B_{1/\varepsilon}} u_\varepsilon^{p-1} v_\varepsilon u_{\varepsilon x}^2 \phi^2 + \int_{B_{1/\varepsilon}} \frac{v_{\varepsilon x}(\cdot,T)^2}{v_\varepsilon (\cdot, T)} \phi^2 +  \int_0^T \int_{B_{1/\varepsilon}}\frac{u_\varepsilon}{v_{\varepsilon}}v_{\varepsilon x}^2 \phi^2 \leq C(p,T).
  $$
\end{proof}
With these results, we can finally bound the sequence $(u_\varepsilon^{\frac{p+1}{2}} v_\varepsilon)_{\varepsilon \in (0,1)}$ in $L^2\big((0,T); W^{1,1}_\text{loc}(\mathbb{R})\big)$ uniformly in $\varepsilon$ for any $T>0$.
\begin{lemma}\label{l4.5}
    Let $p \geq 2$, $K>0$ and assume that $u_0$ and $v_0$ satisfy \eqref{hip0}, \eqref{hip1}, \eqref{3.19} and \eqref{4.1}. Then, for all $T>0$ there exists $C(p,T)>0$ such that 
    $$
    ||u_\varepsilon^{\frac{p+1}{2}} v_\varepsilon||_{ L^2\big((0,T); W^{1,1}_\text{loc}(\mathbb{R})\big)} \leq C(p,T)
    $$
    for all $\varepsilon \in (0,1)$.
\end{lemma}
\begin{proof}
    Given a fixed $T>0$, to prove the local boundedness in space, we restrict the analysis to an arbitrary ball $B_{1/\varepsilon}$. Moreover, we localize the results by means of a cutoff function $\phi^2$ in the sense of Definition \ref{ctf}.
    \\\\
    In this way, for any $\varepsilon \in (0,1)$, we construct a bound for the quantity
    \begin{equation*}
        \begin{split}
        ||u_\varepsilon^{\frac{p+1}{2}} v_\varepsilon \phi^2||&_{ L^2\big((0,T); W^{1,1} (B_{1/\varepsilon}) \big)}^2 = \int_0^T \left(||u_\varepsilon(\cdot,t)^{\frac{p+1}{2}} v_\varepsilon(\cdot,t) \phi^2||_{W^{1,1} (B_{1/\varepsilon})} \right)^2 dt 
          \\\\
          & = \int_0^T\left\{\int_{B_{1/\varepsilon}} u_\varepsilon^{\frac{p+1}{2}} v_\varepsilon \phi^2 + \int_{B_{1/\varepsilon}} \left |\left(u_\varepsilon^{\frac{p+1}{2}} v_\varepsilon \phi^2\right)_x \right| \right\}^2  dt
          \\\\
          & \leq 2 \int_0^T \left\{\int_{B_{1/\varepsilon}} u_\varepsilon^{\frac{p+1}{2}} v_\varepsilon \phi^2 \right\}^2 dt + 2 \int_0^T \left\{\int_{B_{1/\varepsilon}}  \left |\left(u_\varepsilon^{\frac{p+1}{2}} v_\varepsilon \phi^2\right)_x \right|\right\}^2 dt.
        \end{split}
    \end{equation*}
    The first term can be easily bounded for instance by combining Young's inequality with the $L^p$ bound for $u_\varepsilon$ and $L^\infty$ bound for $v_\varepsilon$. In particular, we have
    \begin{equation}\label{4.21}
        \begin{split}
            \left\{\int_{B_{1/\varepsilon}} u_\varepsilon^{\frac{p+1}{2}} v_\varepsilon \phi^2 \right\}^2 & \leq ||v_\varepsilon||_{L^\infty(B_{1/\varepsilon})}^2\left\{\int_{B_{1/\varepsilon}} u_\varepsilon^{\frac{p+1}{2}} \phi^2 \right\}^2 
            \\\\ &\leq \frac{||v_0||_{L^\infty(\mathbb{R})}^2}{4} \left\{\int_{B_{1/\varepsilon}} u_\varepsilon^{p+1} + \int_{B_{1/\varepsilon}} \phi^4 \right\}^2 \leq C_1(p,T),
        \end{split}
    \end{equation}
    for some $C_1(p,T)>0$ provided by Lemma \ref{l3.7}.
    \\\\
    With respect to the second term, one obtains
    \begin{equation*}
            \left\{\int_{B_{1/\varepsilon}} \left |\left(u_\varepsilon^{\frac{p+1}{2}} v_\varepsilon \phi^2\right)_x \right | \right\}^2  \leq 2 \left\{\int_{B_{1/\varepsilon}} \left | \left(u_\varepsilon^{\frac{p+1}{2}} v_\varepsilon \right)_x \right | \phi^2 \right\}^2 +2 \left \{ \int_{B_{1/\varepsilon}} u_\varepsilon^{\frac{p+1}{2}} v_\varepsilon \left |\left( \phi^2 \right)_x \right |\right\}^2,
    \end{equation*}
    where we can estimate the last element as in \eqref{4.21}. For the remaining integral, the Cauchy-Schwarz inequality yields
    \begin{equation}\label{4.22}
        \begin{split} 
        \left\{\int_{B_{1/\varepsilon}} \left | \left(u_\varepsilon^{\frac{p+1}{2}} v_\varepsilon \right)_x \right | \phi^2 \right\}^2 \leq \frac{(p+1)^2}{2} \left\{\int_{B_{1/\varepsilon}} u_\varepsilon^{\frac{p-1}{2}} v_\varepsilon ~ |u_{\varepsilon x}| ~\phi^2 \right\}^2 + 2 \left\{\int_{B_{1/\varepsilon}} u_\varepsilon^{\frac{p+1}{2}} ~ |v_{\varepsilon x}| ~\phi^2 \right\}^2 \\\\
            \leq \frac{(p+1)^2}{2} \left\{\int_{B_{1/\varepsilon}} v_\varepsilon \phi^2 \right\} \cdot \left\{\int_{B_{1/\varepsilon}} u_\varepsilon^{p-1} v_\varepsilon u_{\varepsilon x}^2 \phi^2 \right\} + 2\left\{\int_{B_{1/\varepsilon}} u_\varepsilon^p v_\varepsilon \phi^2 \right\} \cdot \left\{\int_{B_{1/\varepsilon}} \frac{u_\varepsilon}{v_\varepsilon} v_{\varepsilon x}^2 \phi^2 \right\},
        \end{split}
    \end{equation}
    where all the terms are bounded by combining the $L^1$ and $L^\infty$ bounds for $v_\varepsilon$ proved in Lemma \ref{l2.2} with Lemmas \ref{l3.7} and \ref{l4.4}. Therefore, a time integration of the constants allows us to bound $||u_\varepsilon^{\frac{p+1}{2}} v_\varepsilon \phi^2||_{ L^2\big((0,T); W^{1,1} (B_{1/\varepsilon}) \big)}^2$, obtaining boundedness for $||u_\varepsilon^{\frac{p+1}{2}} v_\varepsilon||_{ L^2\big((0,T); W^{1,1}_\text{loc}(\mathbb{R})\big)}$.
\end{proof}
\section{A bound for $\Big(\partial_t(u_\varepsilon^{\frac{p+1}{2}} v_\varepsilon)\Big)_{\varepsilon \in (0,1)}$ in $L^1\big((0,T); (W^{3,2}_\text{loc}(\mathbb{R}))^*\big)$}\label{s5}
In a last step towards establishing a solution $(u,v)$ to the original system \eqref{sist} as a certain limit of the regularized problems, we bound the time derivative $\Big(\partial_t(u_\varepsilon^{\frac{p+1}{2}} v_\varepsilon)\Big)_{\varepsilon \in (0,1)}$ in $L^1\big((0,T); (W^{3,2}_\text{loc}(\mathbb{R}))^*\big)$. In this way, an application of an Aubin-Lions type lemma will allow us to extract a convergent subsequence that defines a the limit solution. We start by proving the following auxiliary result.
\begin{lemma}\label{l5.1}
    Let $K>0$, $q \in (0,1)$, $T>0$ and $\phi^2$ as in Definition \ref{ctf}. Then, if \eqref{hip0}, \eqref{hip1}, and \eqref{3.19} are satisfied, there exists $C(q,T) >0$ such that
    $$
    \int_0^T \int_{B_{1/\varepsilon}} u_\varepsilon^{q-1} v_\varepsilon u_{\varepsilon x}^2 \phi^2 \leq C(q,T),
    $$
    for all $\varepsilon \in (0,1)$.
\end{lemma}
\begin{proof}
    For a fixed $T>0$, integrating by parts we have
    \begin{equation*}
        \begin{split}
            -\frac{1}{q} & \frac{d}{dt} \int_{B_{1/\varepsilon}} u_\varepsilon^q \phi^2 = - \frac{1}{q} \int_{B_{1/\varepsilon}} q u_\varepsilon^{q-1} \left[(u_\varepsilon v_\varepsilon u_{\varepsilon x})_x - (u_\varepsilon^2 v_\varepsilon v_{\varepsilon x})_x + u_\varepsilon v_\varepsilon \right] \phi^2 \\
            & =\int_{B_{1/\varepsilon}} (u_\varepsilon^{q-1} \phi^2)_x (u_\varepsilon v_\varepsilon u_{\varepsilon x}) - \int_{B_{1/\varepsilon}} (u_\varepsilon^{q-1} \phi^2)_x (u_\varepsilon^2 v_\varepsilon v_{\varepsilon x}) - \int_{B_{1/\varepsilon}} u_\varepsilon^q v_\varepsilon \phi^2 \\
            & \leq (q-1) \int_{B_{1/\varepsilon}} u_\varepsilon^{q-1} v_\varepsilon u_{\varepsilon x}^2 \phi^2 + 2\int_{B_{1/\varepsilon}} u_\varepsilon^q v_\varepsilon u_{\varepsilon x} \phi \phi_x \\
            &- (q-1) \int_{B_{1/\varepsilon}} u_{\varepsilon}^q v_\varepsilon u_{\varepsilon x} v_{\varepsilon x} \phi^2 - 2 \int_{B_{1/\varepsilon}} u_\varepsilon^{q+1} v_\varepsilon v_{\varepsilon x} \phi \phi_x, \quad \text{for all } t \in (0,T),
        \end{split}
    \end{equation*}
    where we dropped the non-positive term $- \int_{B_{1/\varepsilon}} u_\varepsilon^q v_\varepsilon \phi^2$. Next, rewriting the above expression we obtain
    \begin{equation}\label{5.1}
        \begin{split}
             -\frac{1}{q} & \frac{d}{dt} \int_{B_{1/\varepsilon}} u_\varepsilon^q \phi^2 + (1-q)\int_{B_{1/\varepsilon}} u_\varepsilon^{q-1} v_\varepsilon u_{\varepsilon x}^2 \phi^2  \leq (1-q) \int_{B_{1/\varepsilon}} u_{\varepsilon}^q v_\varepsilon u_{\varepsilon x} v_{\varepsilon x} \phi^2  \\\\
             & + 2\int_{B_{1/\varepsilon}} u_\varepsilon^q v_\varepsilon u_{\varepsilon x} \phi \phi_x - 2 \int_{B_{1/\varepsilon}} u_\varepsilon^{q+1} v_\varepsilon v_{\varepsilon x} \phi \phi_x, \quad \text{for all } t \in (0,T),
        \end{split}
    \end{equation}
    where we now seek to bound the three terms on the right hand side.
    \\\\
    First, by Young's inequality, we obtain
    \begin{equation}\label{5.2}
        \begin{split}
            (1-q)& \int_{B_{1/\varepsilon}} u_{\varepsilon}^q v_\varepsilon u_{\varepsilon x} v_{\varepsilon x} \phi^2 \leq \frac{1-q}{4} \int_{B_{1/\varepsilon}} u_\varepsilon^{q-1} v_\varepsilon u_{\varepsilon x}^2 \phi^2 + (1-q) \int_{B_{1/\varepsilon}} u_\varepsilon^{q+1} v_{\varepsilon} v_{\varepsilon x}^2 \phi^2 \\\\
            &\leq \frac{1-q}{4} \int_{B_{1/\varepsilon}} u_\varepsilon^{q-1} v_\varepsilon u_{\varepsilon x}^2 \phi^2 + (1-q) ||v_0||_{L^\infty(\mathbb{R})} c_1^2(T) \int_{B_{1/\varepsilon}} u^{q+1}, \quad \text{for all } t \in (0,T),
         \end{split}
    \end{equation}
    where we relied on Lemma \ref{l2.2} for the $L^\infty$ bound for $v_\varepsilon$, Lemma \ref{l3.8} for obtaining a $c_1(T)>0$ such that $||v_{\varepsilon x}||_{L^\infty(B_{1/\varepsilon})} \leq c_1(T)$ and the fact that $\phi^2 \leq 1$ by Definition \ref{ctf}.
    \\\\
    Similarly, for all $t \in (0,T)$ we have
    \begin{equation}\label{5.3}
    \begin{split}
        2&\int_{B_{1/\varepsilon}} u_\varepsilon^q v_\varepsilon u_{\varepsilon x} \phi \phi_x \leq \frac{1-q}{4} \int_{B_{1/\varepsilon}} u_\varepsilon^{q-1} v_\varepsilon u_{\varepsilon x}^2 \phi^2 + \frac{4}{1-q} \int_{B_{1/\varepsilon}} u_\varepsilon^{q+1} v_\varepsilon \phi_x^2 \\\\
        & \leq \frac{1-q}{4} \int_{B_{1/\varepsilon}} u_\varepsilon^{q-1} v_\varepsilon u_{\varepsilon x}^2 \phi^2  + \frac{2 ||v_0||_{L^\infty(\mathbb{R})}}{1-q} \int_{B_{1/\varepsilon}} u_\varepsilon^{2(q+1)} + \frac{2 ||v_0||_{L^\infty(\mathbb{R})}}{1-q} \int_{B_{1/\varepsilon}} \phi_x^4,
    \end{split}
    \end{equation}
    as well as
    \begin{equation}\label{5.4}
        \begin{split}
            - 2 & \int_{B_{1/\varepsilon}} u_\varepsilon^{q+1} v_\varepsilon v_{\varepsilon x} \phi \phi_x \leq 2 \int_{B_{1/\varepsilon}} \left | u_\varepsilon^{q+1} v_\varepsilon v_{\varepsilon x} \phi \phi_x \right | \leq 2 c_1(T)||v_0||_{L^\infty(\mathbb{R})}\int_{B_{1/\varepsilon}} u_\varepsilon^{q+1} \phi |\phi_x|\\
            & \leq c_1(T)||v_0||_{L^\infty(\mathbb{R})}\int_{B_{1/\varepsilon}} u_\varepsilon^{2(q+1)} + c_1(T)||v_0||_{L^\infty(\mathbb{R})}\int_{B_{1/\varepsilon}} \phi^2 \phi_x^2, \quad \text{for all } t \in (0,T).
        \end{split}
    \end{equation}
    Thus, combining the three estimates \eqref{5.2}-\eqref{5.4} and substituting into \eqref{5.1}, we have
    \begin{equation} \label{5.5}
         -\frac{1}{q} \frac{d}{dt} \int_{B_{1/\varepsilon}} u_\varepsilon^q \phi^2 + \frac{(1-q)}{2}\int_{B_{1/\varepsilon}} u_\varepsilon^{q-1} v_\varepsilon u_{\varepsilon x}^2 \phi^2  \leq c_2(q,T) \left (\int_{B_{1/\varepsilon}} u^{q+1} + \int_{B_{1/\varepsilon}} u^{2(q+1)} + 1\right),
    \end{equation}
    for a large enough $c_2(q,T)>0$. Notice that as $1  < q+1 < 2$, a convex interpolation between the $L^1$ and $L^2$ bounds for $u_\varepsilon$ via Young's inequality allows us to obtain $c_3(q,T)>0$ such that $\int_{B_{1/\varepsilon}} u^{q+1} \leq c_3(q,T)$, while the same argument proves the existence of a  $c_4(q,T)>0$ satisfying $\int_{B_{1/\varepsilon}} u^{2(q+1)} \leq c_4(q,T)$, for all $\varepsilon \in (0,1)$, $t \in (0,T)$.
    \\\\
    Thus, a time integration of \eqref{5.5} reveals
    \begin{equation}\label{5.6}
        \frac{1}{q} \int_{B_{1/\varepsilon}} \big(u_0 + \varepsilon \zeta(x)\big)^q \phi^2 + \frac{(1-q)}{2}\int_0^T \int_{B_{1/\varepsilon}} u_\varepsilon^{q-1} v_\varepsilon u_{\varepsilon x}^2 \phi^2 \leq C(q,T) + \frac{1}{q} \int_{B_{1/\varepsilon}} u_\varepsilon^q(\cdot,T) \phi^2,
    \end{equation}
    for all $\varepsilon\in(0,1)$, where $C(q,T):= c_2(q,T) \cdot\Big(c_3(q,T) + c_4(q,T) + 1\Big) \cdot T >0$. 
    \\\\
    The last step is to prove that indeed $\int_{B_{1/\varepsilon}} u_\varepsilon^q(\cdot,T) \phi^2$ is bounded, which can be easily obtained thanks to the presence of the cutoff function. In particular, using Young's inequality with exponents $\frac{q+1}{q}$ and $q+1$ yields
    $$
    \int_{B_{1/\varepsilon}} u_\varepsilon^q(\cdot,T) \phi^2 \leq \frac{q}{q+1} \int_{B_{1/\varepsilon}} u_\varepsilon^{q+1}(\cdot,T) + \frac{1}{q+1}\int_{B_{1/\varepsilon}} \phi^{2(q+1)},    
    $$
    where we have that both quantities are bounded independent of $\varepsilon$, which finishes the proof.
\end{proof}
Once this is proved, we can obtain the desired bound for the time derivative $\Big(\partial_t(u_\varepsilon^{\frac{p+1}{2}} v_\varepsilon)\Big)_{\varepsilon \in (0,1)}$ in $L^1\big((0,T); (W^{3,2}_\text{loc}(\mathbb{R}))^*\big)$.
\begin{lemma}\label{l5.2}
    Let $p \geq 2$, $K>0$ and assume that $u_0$ and $v_0$ are such that \eqref{hip0}, \eqref{hip1}, \eqref{3.19} and \eqref{4.1} hold. Then, for all $T>0$ there exists $C(p,T)>0$ such that
    $$
    \Big|\Big|\partial_t\big(u_\varepsilon^{\frac{p+1}{2}} v_\varepsilon\big)\Big|\Big|_{ L^    1\big((0,T); (W^{3,2}_\text{loc}(\mathbb{R}))^*\big)} \leq C(p,T)
    $$
    for all $\varepsilon \in (0,1)$.
\end{lemma}
\begin{proof}
As in Lemma \ref{l4.5}, given $T>0$, we treat the local boundedness in space by considering an arbitrary ball $B_{1/\varepsilon}$, as well as a cutoff $\phi^2$ as in Definition \ref{ctf}.
\\\\
Thus, for any $\varepsilon \in (0,1)$ we bound the norm of $\partial_t \big(u_\varepsilon^{\frac{p+1}{2}} v_\varepsilon \phi^2 \big)$ in $ L^    1\big((0,T); (W^{3,2}(B_{1/\varepsilon}))^*\big)$ independently of $\varepsilon$. We have
$$
\Big|\Big|\partial_t \big(u_\varepsilon^{\frac{p+1}{2}} v_\varepsilon \phi^2 \big)\Big|\Big|_{ L^1\big((0,T); (W^{3,2}(B_{1/\varepsilon}))^*\big)} = \int_0^T \Big|\Big|\partial_t \big(u_\varepsilon(\cdot, t)^{\frac{p+1}{2}} v_\varepsilon(\cdot, t) \phi^2 \big)\Big|\Big|_{\big(W^{3,2}(B_{1/\varepsilon})\big)^*} ~dt.
$$  
To compute the norm in $\big(W^{3,2}(B_{1/\varepsilon})\big)^*$, we consider an arbitrary $\psi \in W^{3,2}(B_{1/\varepsilon})$ that satisfies $||\psi||_{W^{3,2}(B_{1/\varepsilon})} \leq 1$. In this way, proceed to bound the dual pairing integrating by parts as follows
\begin{equation}\label{5.7}
    \begin{split}
      \int_{B_{1/\varepsilon}} &\partial_t \big(u_\varepsilon^{\frac{p+1}{2}} v_\varepsilon \phi^2 \big) \psi = \frac{p+1}{2} \int_{B_{1/\varepsilon}} u_\varepsilon^{\frac{p-1}{2}} v_\varepsilon \Big((u_\varepsilon v_\varepsilon u_{\varepsilon x})_x - (u_\varepsilon^2 v_\varepsilon v_{\varepsilon x})_x + u_\varepsilon v_\varepsilon \Big) \phi^2 \psi \\\\
      & +\int_{B_{1/\varepsilon}} u_\varepsilon^{\frac{p+1}{2}} \big(v_{\varepsilon xx} - u_\varepsilon v_\varepsilon \big) \phi^2 \psi = - \frac{p+1}{2} \int_{B_{1/\varepsilon}} \Big( u_\varepsilon^{\frac{p-1}{2}} v_\varepsilon \phi^2 \psi\Big)_x (u_\varepsilon v_\varepsilon u_{\varepsilon x} - u_\varepsilon^2 v_\varepsilon v_{\varepsilon x}) \\\\
      &+ \frac{p+1}{2} \int_{B_{1/\varepsilon}} u_\varepsilon^{\frac{p+1}{2}} v_\varepsilon^2 \phi^2 \psi    
      - \int_{B_{1/\varepsilon}} \Big(u_\varepsilon^{\frac{p+1}{2}} \phi^2 \psi \Big)_x v_{\varepsilon x} - \int_{B_{1/\varepsilon}} u_\varepsilon^{\frac{p+3}{2}} v_\varepsilon \phi^2 \psi, \quad \text{for all } t \in (0,T),
    \end{split}
\end{equation}
and computing all the derivatives, one obtains
\begin{equation}\label{5.8}
    \begin{split}
       \int_{B_{1/\varepsilon}} &\partial_t \big(u_\varepsilon^{\frac{p+1}{2}} v_\varepsilon \phi^2 \big) \psi = -\frac{p+1}{2} \cdot \frac{p-1}{2} \int_{B_{1/\varepsilon}} u_\varepsilon^{\frac{p-1}{2}} v_\varepsilon^2 u_{\varepsilon x}^2 \phi^2 \psi - \frac{p+1}{2} \int_{B_{1/\varepsilon}} u_\varepsilon^{\frac{p+1}{2}} v_\varepsilon u_{\varepsilon x} v_{\varepsilon x} \phi^2 \psi \\\\
       & -(p+1)\int_{B_{1/\varepsilon}} u_\varepsilon^{\frac{p+1}{2}} v_\varepsilon^2 u_{\varepsilon x}\phi \phi_x \psi - \frac{p+1}{2} \int_{B_{1/\varepsilon}} u_\varepsilon^{\frac{p+1}{2}} v_\varepsilon^2 u_{\varepsilon x} \phi^2 \psi_x \\\\
       & + \frac{p+1}{2}\cdot \frac{p-1}{2} \int_{B_{1/\varepsilon}} u_\varepsilon^{\frac{p+1}{2}} v_\varepsilon^2 u_{\varepsilon x} v_{\varepsilon x} \phi^2 \psi + \frac{p+1}{2} \int_{B_{1/\varepsilon}} u_\varepsilon^{\frac{p+3}{2}} v_\varepsilon v_{\varepsilon x}^2 \phi^2 \psi \\\\
       & + (p+1) \int_{B_{1/\varepsilon}} u_\varepsilon^{\frac{p+3}{2}} v_\varepsilon^2 v_{\varepsilon x} \phi \phi_x \psi + \frac{p+1}{2} \int_{B_{1/\varepsilon}} u_\varepsilon^{\frac{p+3}{2}} v_\varepsilon^2 v_{\varepsilon x} \phi^2 \psi_x   \\\\
       & - \frac{p+1}{2}\int_{B_{1/\varepsilon}} u_\varepsilon^{\frac{p-1}{2}} u_{\varepsilon x} v_{\varepsilon x} \phi^2 \psi - 2 \int_{B_{1/\varepsilon}} u_\varepsilon^{\frac{p+1}{2}} v_{\varepsilon x} \phi \phi_x \psi - \int_{B_{1/\varepsilon}} u_\varepsilon^{\frac{p+1}{2}} v_{\varepsilon x} \phi^2 \psi_x \\\\
       & + \frac{p+1}{2} \int_{B_{1/\varepsilon}} u_\varepsilon^{\frac{p+1}{2}} v_\varepsilon^2 \phi^2 \psi - \int_{B_{1/\varepsilon}} u_\varepsilon^{\frac{p+3}{2}} v_\varepsilon \phi^2 \psi =: \sum_{i=1}^{13} I_i, \quad \text{for all } t \in (0,T).
    \end{split}
\end{equation}
To conclude the proof, we bound each of these terms individually. It is important to take into account is that, for every fixed $\varepsilon \in (0,1)$, the embedding $W^{3,2}(B_{1/\varepsilon}) \hookrightarrow W^{1, \infty}(B_{1/\varepsilon})$ is continuous, so there exists $c_1(R)>0$ such that $||\psi||_{L^\infty(B_{1/\varepsilon})} + ||\psi_x||_{L^\infty(B_{1/\varepsilon})} \leq c_1(1/\varepsilon)$. However, without extra decay properties, $c_1(1/\varepsilon)$ grows unboundedly as $\varepsilon \to 0$. The key aspect comes from employing the cutoff, $\phi^2$. As both $\phi$ and $\phi_x$ are compactly supported in $B_{1/\varepsilon}$, there does exist $c_2>0$ independent of $\varepsilon$ such that 
\begin{equation} \label{5.9}
    ||\phi^2 \psi||_{L^\infty(B_{1/\varepsilon})} + ||\phi \phi_x \psi||_{L^\infty(B_{1/\varepsilon})} + ||\phi^2 \psi_x||_{L^\infty(B_{1/\varepsilon})} + ||\phi_x \psi ||_{L^\infty(B_{1/\varepsilon})}\leq c_2,
\end{equation}
for all $\varepsilon \in (0,1)$. Moreover, up to a scaling constant, $\phi^2 |\psi|$ behaves essentially as a cutoff function (not necessarily in $C^\infty$ but at least in $C^3$), which makes our previous results applicable to the terms involving it.
\\\\
In this way, starting with $I_1$, we have that for any $q \in (0,1)$
\begin{equation}\label{5.10}
\begin{split}
    |I_1| & \leq \frac{p^2 -1}{4} \int_{B_{1/\varepsilon}} u_\varepsilon^{\frac{p-1}{2}} v_\varepsilon^2 u_{\varepsilon x}^2 \phi^2 |\psi| \leq c(p,q) \cdot \frac{p^2 -1}{4} \int_{B_{1/\varepsilon}} u_\varepsilon^{p-1} v_\varepsilon^2 u_{\varepsilon x}^2 \phi^2 |\psi| \\
    & + c(p,q) \cdot \frac{p^2 -1}{4} \int_{B_{1/\varepsilon}} u_\varepsilon^{q-1} v_\varepsilon^2 u_{\varepsilon x}^2 \phi^2 |\psi| \leq c(p,q) ||v_0||_{L^\infty(\mathbb{R})} \cdot \frac{p^2 -1}{4} \int_{B_{1/\varepsilon}} u_\varepsilon^{p-1} v_\varepsilon u_{\varepsilon x}^2 \phi^2 |\psi| \\
    & + c(p,q) ||v_0||_{L^\infty(\mathbb{R})} \cdot \frac{p^2 -1}{4} \int_{B_{1/\varepsilon}} u_\varepsilon^{q-1} v_\varepsilon u_{\varepsilon x}^2 \phi^2 |\psi|, \quad \text{for all } t \in (0,T),
\end{split}
\end{equation}
where we used that, $q-1 \leq \frac{p-1}{2} \leq p-1$ for any $q \in (0,1)$, soby convexity, Young's inequality provides $c(p,q)>0$ such that $u_\varepsilon^{\frac{p-1}{2}} \leq c(p,q) \Big( u_\varepsilon^{q-1} + u_\varepsilon^{p-1}\Big)$. Notice that, as previously remarked, $\phi^2 |\psi|$ behaves as a cutoff function, and thus the time integral of $|I_1|$ can be bounded by Lemmas \ref{l4.4} and \ref{l5.1}.
\\\\
Next, for $I_2$, first by Lemma \ref{l3.8}, there exists $c_3(T)>0$ with $||v_{\varepsilon x}||_{L^\infty(B_{1/\varepsilon})} \leq c_3(T)$ for all $t \in (0,T)$. Hence, after applying Young's inequality we obtain for all $t \in (0,T)$
\begin{equation}\label{5.11}
    \begin{split}
        |I_2|&  \leq c_3(T) \cdot \frac{p+1}{2} \int_{B_{1/\varepsilon}} u_\varepsilon^{\frac{p+1}{2}} v_\varepsilon |u_{\varepsilon x}| \phi^2 |\psi| \leq  c_3(T) \cdot \frac{p+1}{4} \int_{B_{1/\varepsilon}} u_\varepsilon^{p-1} v_\varepsilon u_{\varepsilon x}^2 \phi^2 |\psi| \\
        & +  c_3(T) \cdot \frac{p+1}{4} \int_{B_{1/\varepsilon}} u_\varepsilon^2 v_\varepsilon \phi^2 |\psi| \leq c_4(p,T) \left(\int_{B_{1/\varepsilon}} u_\varepsilon^{p-1} v_\varepsilon u_{\varepsilon x}^2 \phi^2 |\psi| + 1 \right), 
    \end{split}
\end{equation}
thanks to the $L^\infty$ bound for $v_\varepsilon$ from Lemma \ref{l2.2}, \eqref{5.9} and the $L^2$ bound for $u_\varepsilon$ provided by Lemma \ref{l3.7}.
\\\\
We treat $I_3$ in a similar way, having for all $t \in (0,T)$
\begin{equation}\label{5.12}
\begin{split}
    |I_3| &\leq (p+1)\int_{B_{1/\varepsilon}} u_\varepsilon^{\frac{p+1}{2}} v_\varepsilon^2 |u_{\varepsilon x}|\phi |\phi_x| |\psi| \leq ||v_0||_{L^\infty(\mathbb{R})} \cdot \frac{p+1}{2} \int_{B_{1/\varepsilon}} u_\varepsilon^{p-1} v_\varepsilon u_{\varepsilon x}^2 \phi^2 |\psi| \\
   & + \frac{p+1}{2} ||v_0||_{L^\infty(\mathbb{R})}^2 \cdot  \int_{B_{1/\varepsilon}} u_\varepsilon^2 \phi_x^2 |\psi| \leq c_5(p,T) \left(\int_{B_{1/\varepsilon}} u_\varepsilon^{p-1} v_\varepsilon u_{\varepsilon x}^2 \phi^2 |\psi| + 1  \right).
\end{split}
\end{equation}
For $I_4$ and $I_5$, once again we obtain $c_6(p,T), ~c_7(p,T)>0$ such that for all $t \in (0,T)$
\begin{equation}\label{5.13}
\begin{split}
    |I_4| & \leq \frac{p+1}{2} \int_{B_{1/\varepsilon}} u_\varepsilon^{\frac{p+1}{2}} v_\varepsilon^2 |u_{\varepsilon x}| \phi^2 |\psi_x| \leq ||v_0||_{L^\infty(\mathbb{R})} \cdot \frac{p+1}{4} \int_{B_{1/\varepsilon}} u_\varepsilon^{p-1} v_\varepsilon u_{\varepsilon x}^2 \phi^2 \\
    & +  \frac{p+1}{4} \int_{B_{1/\varepsilon}} u_\varepsilon^2 v_\varepsilon^2 \phi^2 |\psi_x|^2 \leq c_6(p,T) \left(\int_{B_{1/\varepsilon}} u_\varepsilon^{p-1} v_\varepsilon u_{\varepsilon x}^2 \phi^2 + 1  \right). 
    \\\\
    |I_5| & \leq \frac{p^2-1}{4} \int_{B_{1/\varepsilon}} u_\varepsilon^{\frac{p+1}{2}} v_\varepsilon^2 |u_{\varepsilon x}| |v_{\varepsilon x}| \phi^2 |\psi| \leq ||v_0||_{L^\infty(\mathbb{R})} \cdot \frac{p^2-1}{8} \int_{B_{1/\varepsilon}} u_\varepsilon^{p-1} v_\varepsilon u_{\varepsilon x}^2 \phi^2 |\psi| \\
    & + \frac{p^2-1}{8} \int_{B_{1/\varepsilon}} u_\varepsilon^2 v_\varepsilon^2 \phi^2 |\psi| \leq c_7(p,T)  \left(\int_{B_{1/\varepsilon}} u_\varepsilon^{p-1} v_\varepsilon u_{\varepsilon x}^2 \phi^2 |\psi| + 1  \right). 
\end{split}
\end{equation}
For $I_6$, a combination of Young's inequality with the previous $c_3(T)$ bounding $||v_{\varepsilon x}||_{L^\infty(B_{1/\varepsilon})}$ and the $L^p$ bounds for $u_\varepsilon$ yield $c_8(p,T)>0$ such that for all $t \in (0,T)$
\begin{equation}\label{5.14}
    \begin{split}
        |I_6| & \leq \frac{p+1}{2} \int_{B_{1/\varepsilon}} u_\varepsilon^{\frac{p+3}{2}} v_\varepsilon v_{\varepsilon x}^2 \phi^2 |\psi| \leq ||v_0||_{L^\infty(\mathbb{R})}  c_2 c_3(T) \cdot \frac{p+1}{4} \left(\int_{B_{1/\varepsilon}} u_\varepsilon^{p-1} + \int_{B_{1/\varepsilon}} u_\varepsilon^2 \right) \leq c_8(p,T).
    \end{split}
\end{equation}
Next, for $I_7$ and $I_8$ we have $c_9(p,T), ~c_{10}(p,T)>0$ that satisfy for all $t \in (0,T)$
\begin{equation}\label{5.15}
    \begin{split}
        |I_7| &\leq (p+1) \int_{B_{1/\varepsilon}} u_\varepsilon^{\frac{p+3}{2}} v_\varepsilon^2 |v_{\varepsilon x}| \phi |\phi_x |\psi \leq c_2 \cdot c_3(T) \int_{B_{1/\varepsilon}} u_\varepsilon^{\frac{p+3}{2}} v_\varepsilon^2  \\
        & \leq \frac{c_2 c_3(T)||v_0||_{L^\infty(\mathbb{R})}}{2} \left(\int_{B_{1/\varepsilon}} u_\varepsilon^{p-1} + \int_{B_{1/\varepsilon}} u_\varepsilon^4 \right) \leq c_9(p,T).
        \\\\
        |I_8| & \leq \frac{p+1}{2} \int_{B_{1/\varepsilon}} u_\varepsilon^{\frac{p+3}{2}} v_\varepsilon^2 |v_{\varepsilon x}| \phi^2 |\psi_x| \leq c_2 \cdot c_3(T) \cdot \frac{p+1}{2} \int_{B_{1/\varepsilon}} u_\varepsilon^{\frac{p+3}{2}} v_\varepsilon^2 \leq c_{10}(p,T).
    \end{split}
\end{equation}
Concerning $I_9$, again by Young's inequality, and by the estimate provided by Lemma \ref{l4.4}, we have $c_{11}(p,T)>0$ such that
\begin{equation}\label{5.16}
\begin{split}
    |I_9| & \leq  \frac{p+1}{2}\int_{B_{1/\varepsilon}} u_\varepsilon^{\frac{p-1}{2}} |u_{\varepsilon x}| |v_{\varepsilon x}| \phi^2 |\psi| \leq \frac{p+1}{4} \int_{B_{1/\varepsilon}} u_\varepsilon^{p-1} v_\varepsilon u_{\varepsilon x}^2 \phi^2 |\psi| + \frac{p+1}{4} \int_{B_{1/\varepsilon}} \frac{v_{\varepsilon x}^2}{v_\varepsilon} \\
    & \leq c_{11}(p,T) \left(\int_{B_{1/\varepsilon}} u_\varepsilon^{p-1} v_\varepsilon u_{\varepsilon x}^2 \phi^2 |\psi| + 1 \right), \quad \text{for all } t \in (0,T).
\end{split}
\end{equation}
Lastly, from $I_{10}$ to $I_{13}$, the analysis is standard again using that $||\psi||_{W^{3,2}(B_{1/\varepsilon})} \leq 1$, for all $t \in (0,T)$ we have
\begin{equation}\label{5.17}
    \begin{split}
        |I_{10}| & \leq  2 \int_{B_{1/\varepsilon}} u_\varepsilon^{\frac{p+1}{2}} v_{\varepsilon x} \phi \phi_x \psi \leq c_3(T) \int_{B_{1/\varepsilon}} u_\varepsilon^{p+1} \phi^2 \phi_x^2 + \int_{B_{1/\varepsilon}} \phi^2 \leq c_{12}(T),
        \\\\
        |I_{11}| & \leq \int_{B_{1/\varepsilon}} u_\varepsilon^{\frac{p+1}{2}} |v_{\varepsilon x}| \phi^2 |\psi_x | \leq  \frac{c_3(T)} {2} \left( \int_{B_{1/\varepsilon}} u_\varepsilon^{p+1} + \int_{B_{1/\varepsilon}} \phi^4 |\psi_x|^2 \right) \leq c_{13}(T), 
        \\\\
        |I_{12}| & \leq \frac{p+1}{2} \int_{B_{1/\varepsilon}} u_\varepsilon^{\frac{p+1}{2}} v_\varepsilon^2 \phi^2 \psi \leq \frac{p+1}{4} \int_{B_{1/\varepsilon}} u_\varepsilon^{p+1} +||v_0||_{L^\infty(\mathbb{R})}^4\frac{p+1}{4} \int_{B_{1/\varepsilon}} \phi^4 \psi^2 \leq c_{14}(p,T),
        \\\\
        |I_{13}| & \leq \int_{B_{1/\varepsilon}} u_\varepsilon^{\frac{p+3}{2}} v_\varepsilon \phi^2 \psi \leq \frac{1}{2} \int_{B_{1/\varepsilon}} u_\varepsilon^{p+3} + \frac{||v_0||_{L^\infty(\mathbb{R})}^2}{2} \int_{B_{1/\varepsilon}} \phi^4 \psi^2 \leq c_{15}(T),
    \end{split}
\end{equation}
for certain $c_{12}, c_{15} >0$. Integrating these estimates, in view of the properties of Lemmas \ref{l4.4} and \ref{l5.1}, we can conclude that there exists $C(p,T)>0$ verifying
    $$
    \Big|\Big|\partial_t\big(u_\varepsilon^{\frac{p+1}{2}} v_\varepsilon\big)\Big|\Big|_{ L^    1\big((0,T); (W^{3,2}_\text{loc}(\mathbb{R}))^*\big)} \leq C(p,T)
    $$
    for all $\varepsilon \in (0,1)$.
\end{proof}
Lastly, before passing to the limit on our sequence $(u_\varepsilon^{\frac{p+1}{2}} v_\varepsilon)_{\varepsilon\in(0,1)}$, we require one last result that will grant the positivity of $v$, the limit of the sequence $(v_\varepsilon)_{\varepsilon \in (0,1)}$, in order to adequately define the limit $u$.  
\begin{lemma}\label{l5.3}
    Let $K>0$ and assume \eqref{hip0}, \eqref{hip1}, \eqref{3.19} are satisfied. Then, for any $T>0$ there exists $C(T)>0$ such that
    $$
    \int_{B_{1/\varepsilon}} \ln \left(\frac{||v_0||_{L^\infty(B_{1/\varepsilon})}}{v_\varepsilon(\cdot,t)} \right)\phi^2  \leq C(T),
    $$
    for all $t \in (0,T)$, $\varepsilon \in (0,1)$.
\end{lemma}
\begin{proof}
    We compute the derivative of $ \int_{B_{1/\varepsilon}} \ln \left(\frac{||v_0||_{L^\infty(B_{1/\varepsilon})}}{v_\varepsilon(\cdot,t)} \right)\phi^2$ and integrate by parts as follows
    \begin{equation*}
        \begin{split}
            \frac{d}{dt} & \int_{B_{1/\varepsilon}} \ln \left(\frac{||v_0||_{L^\infty(B_{1/\varepsilon})}}{v_\varepsilon} \right)\phi^2 = - \frac{d}{dt} \int_{B_{1/\varepsilon}} \ln v_\varepsilon \phi^2  = - \int_{B_{1/\varepsilon}} \frac{v_{\varepsilon xx}}{v_\varepsilon}\phi^2 + \int_{B_{1/\varepsilon}} u_\varepsilon \phi^2 \\
            & = \int_{B_{1/\varepsilon}} \left(\frac{\phi^2}{v_\varepsilon} \right)_x v_{\varepsilon x} + \int_{B_{1/\varepsilon}} u_\varepsilon \phi^2 = 2 \int_{B_{1/\varepsilon}} \frac{v_{\varepsilon_x}}{v_\varepsilon} \phi \phi_x - \int_{B_{1/\varepsilon}} \frac{v_{\varepsilon x}^2}{v_\varepsilon^2} \phi^2 + \int_{B_{1/\varepsilon}} u_\varepsilon \phi^2, \quad \text{for all } t \in (0,T).
        \end{split}
    \end{equation*}
    Thus, rewriting the expression and using Young's inequality, we obtain
    \begin{equation}\label{5.18}
        \begin{split}
        \frac{d}{dt} & \int_{B_{1/\varepsilon}} \ln \left(\frac{||v_0||_{L^\infty(B_{1/\varepsilon})}}{v_\varepsilon} \right)\phi^2 + \int_{B_{1/\varepsilon}} \frac{v_{\varepsilon x}^2}{v_\varepsilon^2} \phi^2 = 2 \int_{B_{1/\varepsilon}} \frac{v_{\varepsilon_x}}{v_\varepsilon} \phi \phi_x+ \int_{B_{1/\varepsilon}} u_\varepsilon \phi^2 \\
        & \leq \frac{1}{2} \int_{B_{1/\varepsilon}} \frac{v_{\varepsilon x}^2}{v_\varepsilon^2} \phi^2 + 2 \int_{B_{1/\varepsilon}} \phi_x^2 + \int_{B_{1/\varepsilon}} u_\varepsilon \phi^2, \quad \text{for all }t \in (0,T).
        \end{split}
    \end{equation}
    Hence, for all $t \in (0,T)$ one has
    $$
    \frac{d}{dt} \int_{B_{1/\varepsilon}} \ln \left(\frac{||v_0||_{L^\infty(B_{1/\varepsilon})}}{v_\varepsilon} \right)\phi^2 + \frac{1}{2}\int_{B_{1/\varepsilon}} \frac{v_{\varepsilon x}^2}{v_\varepsilon^2} \phi^2 \leq  2 \int_{B_{1/\varepsilon}} \phi_x^2 + \int_{B_{1/\varepsilon}} u_\varepsilon \phi^2\leq c_1,
    $$
    for a given $c_1$ independent of $\varepsilon \in (0,1)$ obtained by the integrability of $u_\varepsilon$ as proved in Lemma \ref{l2.2}. A time integration directly provides the result.
\end{proof}
\section{Passing to the limit $\varepsilon \searrow 0 $} \label{s6}
In this final section, thanks to the estimates obtained for the sequence $(u_\varepsilon^{\frac{p+1}{2}} v_\varepsilon)_{\varepsilon\in(0,1)}$ in Lemmas \ref{l4.5} and \ref{l5.2}, as well as by Lemma \ref{l5.3}, we prove the existence of a subsequence converging to the global weak solution of the original problem \eqref{sist} in the whole real line.
\begin{lemma}\label{l6.1}
    Let $p \geq 1$, assume $u_0$ and $v_0$ satisfy \eqref{hip0} and \eqref{hip01}. Then, there exists a subsequence $(\varepsilon_j)_{j \in \mathbb{N}} \subset (0,1)$ and functions 
    \begin{equation}\label{6.1}
        \begin{cases}
            u \in L^\infty_\text{loc}\big((0,\infty); L^p(\mathbb{R})\big), \\\\
            v \in L^\infty\big((0,\infty); L^p(\mathbb{R})\big) \cap L^\infty\big (\mathbb{R} \times (0, \infty)\big) \cap L^\infty_\text{loc}\big((0,\infty); W^{1,\infty}(\mathbb{R})\big), 
        \end{cases}
    \end{equation}
    such that $u\geq 0$ in $\mathbb{R} \times (0,\infty)$, $v > 0$ in $\mathbb{R} \times (0,\infty)$ and $\varepsilon_j \searrow 0$ as $j \to \infty$ satisfying
    \begin{align}
        & u_\varepsilon \to  u  \quad \text{a.e. in } \mathbb{R} \times (0,\infty) \text{ and in } L^q_\text{loc} (\mathbb{R} \times (0,\infty)) \text{ for all } q \in [1,p),  \label{6.2} \\
        & v_\varepsilon \to v \quad \text{a.e. in } \mathbb{R} \times (0,\infty) \text{ and in } L^q_\text{loc} (\mathbb{R} \times (0,\infty)) \text{ for all } q \in [1,\infty), \label{6.3} \\
        & v_{\varepsilon x} \overset{\ast}{\rightharpoonup}  v_x  \quad \text{ in } L^\infty_{\text{loc}}\big((0,\infty); L^\infty(\mathbb{R)\big)}, \label{6.4}
    \end{align}
    as $\varepsilon = \varepsilon_j \searrow 0$. Moreover, the pair $(u,v)$ is a global weak solution to system \eqref{sist} in the sense of Definition \ref{weak-sol}.
\end{lemma}
\begin{proof}
    We start by the convergence of $v_\varepsilon$. By Lemma \ref{l2.2} we know that both $||v_\varepsilon||_{L^1(B_{1/\varepsilon})}$ and $||v_\varepsilon||_{L^\infty(B_{1/\varepsilon})}$ are uniformly bounded for all $t >0$ and $\varepsilon \in (0,1)$. Moreover, by Lemma \ref{l3.8} for any $T>0$, $||v_{\varepsilon x}||_{L^\infty(B_{1/\varepsilon})}$ is also uniformly bounded in $\varepsilon$ for all $t \in (0,T)$.
    \\\\
    Hence, this provides bounds for $(v_\varepsilon)_{\varepsilon \in (0,1)}$ in $L^\infty \big((0, \infty); L^p (B_{1/\varepsilon}) \big)$ for any given $p \geq 1$, as well as in $L^\infty\big(B_{1/\varepsilon} \times (0,\infty) \big)$ and in $L^\infty \big((0, T); W^{1,\infty} (B_{1/\varepsilon}) \big)$ uniformly in $\varepsilon \in (0,1)$. In particular, this implies that the sequence $(v_\varepsilon)_{\varepsilon \in (0,1)}$ is bounded in $L^q\big((0, T); W^{1,q}_{\text{loc}} (\mathbb{R}) \big)$ for any $q \in [1, \infty)$ and any $T>0$ uniformly in $\varepsilon$.
    \\\\
    Furthermore, it is direct to check that $(v_{\varepsilon t})_{\varepsilon \in (0,1)}$ is bounded in $L^2\big((0,T); (W^{1,2}_\text{loc}(\mathbb{R}))^*\big)$. This can be proved easily by taking into account that the non-linear term $-u_\varepsilon v_\varepsilon$ appearing in $v_{\varepsilon t}$ is bounded in $L^2_\text{loc}\big(\mathbb{R} \times (0,\infty) \big)$, as well as a standard duality argument for $v_{\varepsilon xx}$.
    \\\\
    Thus, by means of an Aubin-Lions type lemma \cite{aubin-lions}, we can obtain $(\varepsilon_{j_1})_{j_1 \in \mathbb{N}} \subset (0,1)$ with $\varepsilon_{j_1} \searrow 0$ as $j_1 \to \infty$ such that there exists 
    $$
    v \in L^\infty\big((0,\infty); L^p(\mathbb{R})\big) \cap L^\infty\big (\mathbb{R} \times (0, \infty)\big) \cap L^\infty_\text{loc}\big((0,\infty); W^{1,\infty}(\mathbb{R})\big),  \quad \text{for all } p \geq 1,
    $$
    satisfying $v_\varepsilon \to v$ for $\varepsilon = \varepsilon_{j_1}$ almost everywhere in $\mathbb{R} \times (0, \infty)$ and in $L^q_\text{loc} (\mathbb{R} \times (0,\infty)) \text{ for all } q \in [1,\infty)$. Moreover, the strict positivity of $v$ follows from combining the fact that $v \leq ||v_0||_{L^\infty(\mathbb{R})}$, Lemma \ref{l5.3} and Fatou's lemma, ensuring that $\ln v \in L^1_{\text{loc}}\big(\mathbb{R} \times (0,\infty)\big)$ and indeed $v$ is positive almost everywhere in $\mathbb{R} \times (0,\infty)$.
    \\\\
    The weak$^*$ convergence of $v_{\varepsilon x}$ follows from the Banach-Alaoglu theorem, as a consequence of the boundedness of $||v_{\varepsilon x}||_{L^\infty(B_{1/\varepsilon})}$ for all $t \in (0,T)$ for arbitrary $T>0$. The strong convergence of $v_\varepsilon$ guarantees that indeed the limit of $v_{\varepsilon x}$ coincides with $v_x$, granting \eqref{6.4}.
    \\\\
    Next, we assess convergence of $u_\varepsilon$ by means of the results derived for the sequence $(u_\varepsilon^\frac{p+1}{2} v_\varepsilon)_{\varepsilon \in (0,1)}$ gathered in Lemma \ref{l4.5} and Lemma \ref{l5.2}. Specifically, for arbitrary $p \geq 2$ and $T > 0$, we have
    \begin{equation*}
    \begin{split}
        (u_\varepsilon^\frac{p+1}{2} v_\varepsilon)_{\varepsilon \in (0,1)} \quad & \text{ is bounded in } L^2\big((0,T); W^{1,1}_{\text{loc}}(\mathbb{R})\big), \\
     \Big(\partial_t(u_\varepsilon^{\frac{p+1}{2}} v_\varepsilon)\Big)_{\varepsilon \in (0,1)} \quad &  \text{ is bounded in } L^1\big((0,T); (W^{3,2}_\text{loc}(\mathbb{R}))^*\big).
    \end{split}
    \end{equation*}
    In particular, these bounds hold for the subsequence $\varepsilon = \varepsilon_{j_1}$ previously extracted to ensure the convergence of $v_\varepsilon$. Thus, a further application of the Aubin-Lions lemma provides a second subsequence $(\varepsilon_{j_2})_{j_2 \in \mathbb{N}} \subset (\varepsilon_{j_1})_{j_1 \in \mathbb{N}} \subset (0,1)$ with $\varepsilon_{j_2} \searrow 0$ as $j_2 \to \infty$ such that there exists 
    $$
    z \in L^1_{\text{loc}}\big(\mathbb{R} \times (0, \infty)\big),  ~ \text{ satisfying } ~ u_\varepsilon^\frac{p+1}{2} v_\varepsilon \to z ~\text{ a. e. in } \mathbb{R} \times (0,\infty) \text{ and in } L^1_{\text{loc}}\big(\mathbb{R} \times (0, \infty)\big),
    $$
    for $\varepsilon = \varepsilon_{j_2}$. In this way, thanks to the positivity of $v$, we can define the limit $u :=\displaystyle \left( \frac{z}{v} \right)^{\frac{2}{p+1}}$. Moreover, by Lemma \eqref{2.2} and Lemma \ref{l3.7}, $(u_\varepsilon)_{\varepsilon \in (0,1)}$ is bounded in $L^\infty_{\text{loc}}\big( (0, \infty); L^p(B_{1/\varepsilon}) \big)$ for all $p \geq 1$ uniformly in $\varepsilon \in (0,1)$, which together with the Vitali convergence theorem implies that $u \in  L^\infty_\text{loc}\big((0,\infty); L^p(\mathbb{R})\big)$ with
    $u_\varepsilon \to u$ almost everywhere in  $ \mathbb{R} \times (0,\infty)$ and in  $L^q_{\text{loc}} \big(\mathbb{R} \times (0,\infty) \big)$ for all $q \in [1,p)$.
    \\\\
    The last step is to check that $(u,v)$ indeed solve system \eqref{sist} in the sense specified in Definition \ref{weak-sol}. To this end, we consider a test function $\varphi \in C_0^\infty (\mathbb{R} \times [0, \infty))$. As for all $\varepsilon \in (0,1)$, $(u_\varepsilon, v_\varepsilon)$ classically solve the regularized system \eqref{reg}, in particular we have
    \begin{equation*}
    \begin{split}
    \int_0^\infty \int_{B_{1/\varepsilon}} u_\varepsilon \varphi_t + \int_{B_{1/\varepsilon}} u_0 \varphi(\cdot, 0) = &  -\frac{1}{2}\int_0^\infty \int_{B_{1/\varepsilon}} u_\varepsilon^2 v_{\varepsilon x} \varphi_x - \frac{1}{2}\int_0^\infty \int_{B_{1/\varepsilon}} u_\varepsilon v_\varepsilon \varphi_{xx} \\
    & - \int_0^\infty \int_{B_{1/\varepsilon}} u_\varepsilon^2 v_\varepsilon v_{\varepsilon x} \varphi_x  - \int_0^\infty \int_{B_{1/\varepsilon}} u_\varepsilon v_\varepsilon \varphi,
    \end{split}
    \end{equation*}
and
$$
    \int_0^\infty \int_{B_{1/\varepsilon}} v_\varepsilon \varphi_t + \int_{B_{1/\varepsilon}} v_0 \varphi(\cdot, 0) =\int_0^\infty \int_{B_{1/\varepsilon}} v_{\varepsilon x} \varphi_x + \int_0^\infty \int_{B_{1/\varepsilon}} u_\varepsilon v_\varepsilon \varphi,
$$
Therefore, by passing to the limit with $\varepsilon = \varepsilon_{j_2} \searrow 0$, given the convergence properties \eqref{6.2}-\eqref{6.4}, we deduce that the limit functions $(u,v)$ satisfy the weak formulation of system \eqref{sist} in the sense of Definition \ref{weak-sol}, which concludes the proof.
\end{proof}
This allows us to lastly prove the main result.
\\\\
Proof of Theorem \ref{t1}. We only need to employ Lemma \ref{l5.3}. Notice that hypotheses \eqref{hip1}, \eqref{3.19} and \eqref{4.1} assumed for the previous lemmas are a direct consequence of \eqref{hip0} and \eqref{hip01}.
\section*{Acknowledgements}
This work was partially supported by Grant FPU23/03170 from the Spanish Ministry of Science, Innovation and Universities, as well as by a DAAD Research Grant - Bi-nationally Supervised Doctoral Degrees/Cotutelle, 2024/25 (57693451), Grant number 91907995 (F.H.H.).

\end{document}